\documentclass[11pt]{amsart}
\usepackage[margin=1in]{geometry}
\usepackage{amsmath,amssymb,mathtools}
\usepackage{nicefrac}
\usepackage[square]{natbib}
\usepackage[dvipsnames]{xcolor}
\usepackage{xurl}
\newcommand{\bfa}{\mathbf a}
\usepackage{tikz}
\usetikzlibrary{arrows.meta,shapes.misc,decorations.pathreplacing}
\usepackage[citecolor=PineGreen,colorlinks=true,linkcolor=RoyalBlue,urlcolor=BrickRed]{hyperref}

\newcommand{\nf}{\nicefrac}
\newcommand{\Z}{\mathbb{Z}}
\newcommand{\R}{\mathbb{R}}
\newcommand{\C}{\mathbb{C}}
\newcommand{\T}{\mathbb{T}}
\newcommand{\PP}{\mathbb{P}}
\newcommand{\E}{\mathbb{E}}
\DeclareMathOperator{\adj}{adj}
\newcommand{\pgram}[1]{\mathord{\tikz[baseline=0pt]{\node[inner sep=0pt,outer sep=0pt,anchor=base,text height=1ex,text depth=0.2em](p) at (0,0){$#1$};\draw[line width=0.4pt,line join=miter]([xshift=-0.36em,yshift=-0.02em]p.south west)--([xshift=0.08em,yshift=-0.02em]p.south east)--([xshift=0.36em,yshift=0.24em]p.north east)--([xshift=-0.08em,yshift=0.24em]p.north west)--cycle;}}}

\newtheorem{theorem}{Theorem}[section]

\newtheorem{lemma}[theorem]{Lemma}
\newtheorem{corollary}[theorem]{Corollary}
\theoremstyle{remark}

\numberwithin{equation}{section}

\title[Stationary harmonic functions]{Stationary harmonic functions on lattices are constant}
\author{Ahmed Bou-Rabee}
\date{}

\begin{document}

\begin{abstract}
We prove that a random harmonic function on $\mathbb Z^d$ whose law is invariant under translations is constant. The proof reduces the problem to a deterministic estimate: a harmonic function that crosses the interval $(-1,1)$ along a positive fraction of the edges of a box is exponentially large on a positive fraction of a larger box. Both results extend to periodic networks. For uniformly elliptic conductances, the theorem holds in dimension two but fails in higher dimensions, even when the conductances and the function are jointly stationary.
\end{abstract}

\maketitle

\section{Introduction}\label{s.intro}

Let $d\geq1$ and let $e_1,\dots,e_d$ be the standard basis vectors of $\Z^d$. A function $u\colon\Z^d\to\R$ is \textbf{harmonic} if $\Delta u=0$, where
\begin{equation*}
\Delta u(x)\coloneqq\frac1{2d}\sum_{i=1}^d\bigl(u(x+e_i)+u(x-e_i)-2u(x)\bigr)\,.
\end{equation*}
A \textbf{random harmonic function} is a random function $h\colon\Z^d\to\R$ with $\Delta h=0$ almost surely. We say that $h$ is \textbf{stationary} if, for every $y\in\Z^d$, the translate $(h(x+y))_{x\in\Z^d}$ has the same law as $h$.

In their work on the divisible sandpile, \citet*{LMPU} asked when a random harmonic function is almost surely constant. The answer depends on the graph. Every Cayley graph with nontrivial Poisson boundary carries a nonconstant stationary harmonic function \citep*[Theorem~2]{HLM}. In contrast, the main theorems of \citet*{BLMS} and \citet*{BCG} imply that stationary harmonic functions on $\Z^2$, and more generally on periodic planar graphs with periodic conductances, are constant. It was unclear whether the same holds on nonplanar graphs, in particular on $\Z^d$ for $d\geq3$.\footnote{The question is related to the Gaussian free field: \citet*[Section~2.4]{AGP} recently noted that, for $d\geq3$, the Gaussian free field on $\Z^d$ has an ergodic Gibbs measure with infinite specific free energy if and only if there is a nonconstant random harmonic function on $\Z^d$ whose law is ergodic under translations.}

\subsection{Main results}\label{s.results}

Our main result answers this question for $\Z^d$ in every dimension.

\begin{theorem}\label{t.Zd}
Let $d\geq1$ and let $h$ be a stationary random harmonic function on $\Z^d$. Then almost surely $h(x)=h(0)$ for every $x\in\Z^d$.
\end{theorem}

We deduce Theorem~\ref{t.Zd} from a deterministic estimate: a harmonic function that crosses the interval $(-1,1)$ along a positive fraction of the edges of a box is exponentially large on a positive fraction of a larger box. For real $n\geq0$, let $Q_n\coloneqq[-n,n]^d\cap\Z^d$. We say that a property holds for at least a fraction $p$ of the elements of a finite set $A$ if it holds for at least $p|A|$ of them.

\begin{theorem}\label{t.growth}
Let $d\geq1$. There are $\rho>0$ and $b>0$, depending only on $d$, such that for every $p\in(0,1]$ there are $\kappa>0$ and $N_0$, depending only on $d$ and $p$, with the following property. Let $N\geq N_0$, let $u\colon\Z^d\to\R$ be harmonic, and suppose that at least a fraction $p$ of the edges $\{x,y\}\subset Q_{\rho N}$ satisfy
\begin{equation}\label{e.pair}
u(x)\,u(y)<0\qquad\text{and}\qquad|u(x)|,\,|u(y)|\geq1\,.
\end{equation}
Then $|u(x)|>e^{bpN}$ for at least a fraction $\kappa$ of the sites $x\in Q_{3N}$.
\end{theorem}

\begin{proof}[Proof of Theorem~\ref{t.Zd}, assuming Theorem~\ref{t.growth}]
Suppose that $h$ is nonconstant with positive probability. By stationarity and a normalization, we may assume that~\eqref{e.pair} holds with positive probability on some edge $\{0,e_j\}$, so that, by stationarity, the expected fraction of the edges of $Q_{\rho N}$ that satisfy~\eqref{e.pair} is a constant $2p>0$. Hence, with probability at least $p$, at least a fraction $p$ of these edges satisfy~\eqref{e.pair}, and then Theorem~\ref{t.growth} gives $|h|>e^{bpN}$ on at least a fraction $\kappa$ of $Q_{3N}$, for $N\geq N_0$. Taking expectations and using stationarity, we get $\PP(|h(0)|>e^{bpN})\geq p\kappa$ for $N\geq N_0$, which contradicts $\PP(|h(0)|>e^{bpN})\to0$.
\end{proof}

The proof uses stationarity only through the fact that, for each $j$, the law of $(h(x),h(x+e_j))$ does not depend on $x$, and it works with every lower bound in Theorem~\ref{t.growth} that tends to infinity with $N$.

If $\E|h(0)|<\infty$, then Theorem~\ref{t.Zd} follows from an exercise on conditional expectations \citep[see, for example,][Exercise~4.1.10]{Durrett}, and \citet[Lemma~5.4]{LMPU} proved it when $\E h(0)^+<\infty$.

The closest result to Theorem~\ref{t.growth} is a lower bound of~\citet{BLMS} on $\Z^2$: if a harmonic function $u$ satisfies $|u|\geq2$ somewhere near the origin and $|u|\leq1$ on most sites at every scale, then $|u|$ is exponentially large somewhere \citep[Theorem~(B)]{BLMS}. \citet[Theorem~A.1]{BCG} proved the same lower bound on periodic planar graphs with periodic conductances. Theorem~\ref{t.growth} differs from this lower bound in both its hypothesis and its conclusion. It assumes a sign change across a positive fraction of the edges of a single box, with $|u|\geq1$ at both endpoints, rather than $|u|\leq1$ on most sites at every scale. It concludes that $|u|$ is exponentially large on a positive fraction of a larger box, rather than at one site. \citet{BLMS} combined their lower bound with an upper bound to prove a Liouville theorem: for some $\varepsilon>0$, a harmonic function $u$ on $\Z^2$ is constant if $|u|\leq1$ on all but a fraction $\varepsilon$ of the sites of $Q_N$ for every large $N$ \citep[Theorem~1.1]{BLMS}. \citet[Theorem~1.1]{BCG} proved a Liouville theorem of this kind on periodic planar graphs with periodic conductances, using a different, more geometric argument for the upper bound.

The Liouville theorem of \citet{BLMS} fails on $\Z^d$ for $d\geq3$. For example, as they noted, the function $(3+2\sqrt2)^{x_3}(-1)^{x_1}\mathbf 1\{x_1=x_2\}$ is harmonic on $\Z^3$ and vanishes off the plane $\{x_1=x_2\}$. In contrast, Theorem~\ref{t.growth} holds in every dimension.

In the continuum, there are nonconstant stationary harmonic functions on $\R^2$ \citep{Weiss,Glucksam,BGLS,GlucksamWeiss}. On $\Z^2$, \citet{Tsirelson} constructed a stationary random vector field whose divergence equals $1$ almost surely.

\subsection{Periodic networks}\label{s.periodic-results}

Both theorems extend to periodic networks. Let $G$ be a connected, locally finite graph with vertex set $V$ and positive conductances on its edges; we allow loops and multiple edges. Let $\bfa(v,w)$ be the total conductance of the edges between $v$ and $w$, and let $\bfa(v)\coloneqq\sum_{w\in V}\bfa(v,w)$. A function $u\colon V\to\R$ is harmonic if $\Delta u=0$, where $\Delta u(v)\coloneqq \bfa(v)^{-1}\sum_{w\in V}\bfa(v,w)(u(w)-u(v))$. Suppose that $\Z^d$ acts freely on $V$ by bijections $\tau_x$ that preserve the conductances, and that the action has finitely many orbits. Examples include the triangular and hexagonal lattices, and $\Z^d$ with conductances that are periodic in each coordinate. In place of $u$, we work with the function $x\mapsto u(\tau_xo)$ on $\Z^d$ for a fixed vertex $o$, so that $Q_n$ and the edges $\{x,y\}$ are those of $\Z^d$.

\begin{theorem}\label{t.periodic}
Let $G$, the conductances and the action be as above, and fix $o\in V$. There are $\rho>0$ and $b>0$, depending only on $G$, the conductances, the action and $o$, such that for every $p\in(0,1]$ there are $\kappa>0$ and $N_0$, depending also on $p$, with the following property. Let $N\geq N_0$, let $u\colon V\to\R$ be harmonic, and suppose that at least a fraction $p$ of the edges $\{x,y\}\subset Q_{\rho N}$ satisfy~\eqref{e.pair} for the function $x\mapsto u(\tau_xo)$. Then $|u(\tau_xo)|>e^{bpN}$ for at least a fraction $\kappa$ of the sites $x\in Q_{3N}$.
\end{theorem}

A random harmonic function on $G$ is stationary if its law is invariant under the action.

\begin{corollary}\label{c.periodic}
Let $G$, the conductances and the action be as above, and let $h$ be a stationary random harmonic function on $G$. Then almost surely $h(v)=h(w)$ for all $v,w\in V$.
\end{corollary}

\begin{proof}[Proof of Corollary~\ref{c.periodic}, assuming Theorem~\ref{t.periodic}]
For each $o\in V$, the random function $x\mapsto h(\tau_xo)$ on $\Z^d$ is stationary, so the proof of Theorem~\ref{t.Zd}, with Theorem~\ref{t.periodic} in place of Theorem~\ref{t.growth}, shows that it is almost surely constant. Then $h$ induces a harmonic function on the finite connected quotient graph $G/\Z^d$, with the induced conductances, so it is constant \citep[see, for example,][Lemma~1.16]{LPW}.
\end{proof}

\subsection{Uniformly elliptic conductances}\label{s.random-results}

For uniformly elliptic conductances, the answer depends on the dimension. Let $\bfa(x,y)=\bfa(y,x)>0$ be conductances on the nearest-neighbor edges of $\Z^d$. A function $u\colon\Z^d\to\R$ is \textbf{$\bfa$-harmonic} if $\sum_{|y-x|=1}\bfa(x,y)(u(y)-u(x))=0$ for every $x$.

\begin{theorem}\label{t.random}\leavevmode
\begin{enumerate}
\item[(a)] Let $\lambda\geq1$, and let $h$ be a stationary random function on $\Z^2$ that is almost surely $\bfa$-harmonic for some conductances $\bfa$ with $\lambda^{-1}\leq\bfa\leq\lambda$ that depend measurably on the randomness. Then almost surely $h(x)=h(0)$ for every $x\in\Z^2$.
\item[(b)] There are random conductances $\bfa\leq1$ and a random function $h$ on $\Z^2$ such that $(\bfa,h)$ is jointly stationary and $h$ is almost surely $\bfa$-harmonic and not constant. There are also such $\bfa$ and $h$ with $\bfa\geq1$.
\end{enumerate}
\end{theorem}

In dimensions three and higher, Theorem~\ref{t.random}(a) fails, even when $(\bfa,h)$ is jointly stationary.

\begin{theorem}\label{t.counter}
Let $d\geq3$ and $\lambda>1$. There are random conductances $\bfa$ and a random function $h$ on $\Z^d$ such that $(\bfa,h)$ is jointly stationary, $\lambda^{-1}<\bfa<\lambda$, and $h$ is almost surely $\bfa$-harmonic and not constant.
\end{theorem}

\subsection{Outline of the proof}\label{s.idea}

The proof of Theorem~\ref{t.growth} starts with an observation: its hypothesis cannot be satisfied by a polynomial of degree less than a small multiple of $pN$. Indeed, some coordinate line would then contain more edges across which the polynomial changes sign than its degree, and the restriction of the polynomial to this line has a root inside each of these edges.

To pass from polynomials to harmonic functions, we first use large-scale analyticity, introduced by~\citet*{AKS} for periodic elliptic equations. For a harmonic function $u$ and every integer $1\leq m\leq cN$, it gives a polynomial $P$ of degree at most $m$ such that
\begin{equation*}
\max_{Q_{cN}}|u-P|\leq2^{-m}\max_{Q_{3N}}|u|\,,
\end{equation*}
where $c>0$ depends only on $d$. \citet{BCG} proved this estimate on periodic planar graphs, and in Appendix~\ref{s.analytic} we explain how their proof applies to $\Z^d$. Now let $u$ be as in Theorem~\ref{t.growth}, and let $m$ be a small multiple of $pN$. If $\max_{Q_{3N}}|u|<2^m$, then $P$ changes sign across every edge where~\eqref{e.pair} holds, and this is impossible for a polynomial of degree at most $m$. Hence $\max_{Q_{3N}}|u|\geq2^m$.

The main new step, needed for the proof of Theorem~\ref{t.Zd}, is to show that $|u|$ is exponentially large on a positive fraction of $Q_{3N}$. We argue by contradiction. If $|u|$ is exponentially large at only a few sites, then a nonzero polynomial $F$ of small degree vanishes at all of them. The product $Fu$ is then not exponentially large, and $\Delta^r(Fu)=0$ for some $r$. In Section~\ref{s.products}, we show that such functions are also close to polynomials, by extending them to harmonic functions in one more dimension. So $Fu$ is close to a polynomial $P$, and $FP$ has the sign of $u$ at every site where $|u|\geq1$ and $|F|$ is not small. In Section~\ref{s.proof}, we use the Remez inequality to show that $|F|$ is small compared with its maximum at only a few sites. Hence $FP$ changes sign across more edges than its degree, a contradiction.

Theorem~\ref{t.periodic} follows from the same argument, once we reduce harmonicity on a periodic network to an equation with constant coefficients on each orbit (Section~\ref{s.periodic}). For uniformly elliptic conductances, no such reduction is available, and we use different arguments, outlined at the beginning of each part of Section~\ref{s.random}. In Section~\ref{s.planar}, we prove Theorem~\ref{t.random}(a) by planarity and elliptic regularity. In Sections~\ref{s.examples} and~\ref{s.higher}, we construct the counterexamples of Theorem~\ref{t.random}(b) and Theorem~\ref{t.counter}.

\subsection{Notation and conventions}\label{s.notation}

\begin{itemize}
\item We write $|x|$, $|x|_1$ and $|x|_\infty$ for the Euclidean, $\ell^1$ and maximum norms of $x\in\Z^d$.
\item We write $Q_n^d$ for $Q_n$ when we need to indicate the dimension.
\item The boundary $\partial A$ of a set $A\subseteq\Z^d$ is the set of sites outside $A$ that have a neighbor in $A$.
\item For a finite set $Q$ and a real-valued function $f$ defined on a set containing $Q$, we write
\begin{equation*}
\|f\|_{L^\infty(Q)}\coloneqq\max_Q|f|\,,\qquad\|f\|_{L^2(Q)}\coloneqq\biggl(\sum_{x\in Q}|f(x)|^2\biggr)^{\nf12}\qquad\text{and}\qquad\|f\|_{\underline L^2(Q)}\coloneqq|Q|^{-\nf12}\|f\|_{L^2(Q)}\,.
\end{equation*}
\item The constant $C$ is positive and may change from line to line. We indicate its dependence by writing, for example, $C(d)$.
\end{itemize}

\subsection{Formalization}\label{s.formalization}

A Lean~4 formalization of this paper is maintained at
\begin{center}
\url{https://github.com/nitromannitol/stationary-harmonic-lean}
\end{center}
It verifies every theorem, corollary and lemma of this paper. It also proves the following extensions, among others listed in its documentation:
\begin{itemize}
\item Theorem~\ref{t.random}(a) holds on every connected, locally finite planar graph that is invariant under a lattice of translations isomorphic to $\Z^2$, with finitely many vertex orbits.
\item Conversely, Theorem~\ref{t.random}(a) fails, for some $\lambda$, on every nonplanar Cayley graph of $\Z^2$ with a finite generating set.
\item A difference equation with constant coefficients on $\Z^d$ has only constant stationary solutions if and only if its symbol has no zero on the unit torus other than $(1,\dots,1)$. Theorems~\ref{t.Zd} and~\ref{t.growth} also hold for powers of the Laplacian.
\item Theorem~\ref{t.random}(a) holds in every dimension when the conductances are random and independent of $h$, or deterministic and not necessarily periodic.
\end{itemize}

\subsection*{Use of artificial intelligence}

The author has worked on this problem without success since 2019 \citep[see, for example, the related results in][]{BCD,BCG}, and tried it on every publicly available large language model. When Gemini~4 Carbon was released, the author noticed that it had good ideas but struggled to complete proofs. The author therefore set up a loop between models: Gemini~4 Carbon agents proposed proofs and ideas, and GPT-6.2 Astra agents wrote complete arguments from them, checked every step and reported gaps back through a shared channel. The results of this paper were developed in this loop.

The first complete proof found in this loop was for $\Z^3$. It was long and specific to three dimensions, and it used the harmonic functions supported on diagonal planes. The author then asked the agents to extend this proof to all dimensions and to periodic graphs. They found a shorter but less transparent proof that covers every periodic graph in every dimension, based on Fourier analysis and Siegel's lemma. After reading these arguments, the author proposed that the agents explore an alternative proof: to combine large-scale analyticity with the discrete Remez inequality, as in the appendix of~\citet{BCG}, and to reduce the problem to a deterministic estimate. The author then worked through the details of these arguments with the models to find the proof of Theorems~\ref{t.growth} and~\ref{t.periodic} presented in this paper.

After this, the author asked the loop whether the situation changes for uniformly elliptic conductances, giving it some earlier unpublished partial work of the author. It proved the two-dimensional case and found the counterexamples in dimensions three and higher (Section~\ref{s.higher}). Its proof of the two-dimensional case was not very transparent and used the first-order corrector of stochastic homogenization. The author then asked the agents for a more direct proof of the two-dimensional case, using level sets and planarity as in~\citet{BCG}, which led to the proof in Section~\ref{s.planar}.

The exposition was written by the author; Claude assisted with typesetting and with wording. The author takes responsibility for the mathematical content and presentation of the paper.

\subsection*{Acknowledgments}

Thank you to Jeanne Boursier, William Cooperman, Shirshendu Ganguly, Lionel Levine, Alexander Logunov, Christoforos Panagiotis, Ron Peled, Charles Smart, Mikhail Sodin, and Lingfu Zhang for useful conversations over the years. Thank you to Alexander Holroyd, Lionel Levine, and Mathav Murugan for sharing their unpublished manuscript on Cauchy harmonic functions.

\section{Large-scale analyticity for polyharmonic functions}\label{s.products}

We first recall the Remez inequality.

\begin{lemma}[{\citealp[Lemma~2.1]{BLMS}}]\label{l.remezBLMS}
Let $F\colon\R\to\R$ be a real polynomial of degree at most $k$, let $J$ be a closed interval, and let $E\subset J$ be measurable with positive length $|E|$. Then
\begin{equation*}
\max_J|F|\leq\biggl(\frac{4|J|}{|E|}\biggr)^k\sup_E|F|\,.
\end{equation*}
\end{lemma}

We extend large-scale analyticity for harmonic functions, Lemma~\ref{l.harmonic}, to polyharmonic functions, that is, functions $v$ with $\Delta^rv=0$ for some integer $r\geq1$. These include the products of harmonic functions with polynomials, which we use in Section~\ref{s.proof}.

\begin{lemma}\label{l.analyticity}
There are $\rho>0$ and $N_0$, depending only on $d$, with the following property. Let $N\geq N_0$, let $m$ and $r$ be integers with $1\leq m\leq\rho N$ and $1\leq r\leq m+1$, and let $v\colon\Z^d\to\R$ satisfy $\Delta^rv=0$ in $Q_{3N}$. There is a real polynomial $P\colon\R^d\to\R$ of total degree at most $m$ such that
\begin{equation*}
\|v-P\|_{L^\infty(Q_{\rho N})}\leq2^{-m}\|v\|_{L^\infty(Q_{3N})}\,.
\end{equation*}
\end{lemma}

\begin{proof}
We extend $v$ to a harmonic function $W$ in one more dimension, bound $W$ in $L^2$, and apply large-scale analyticity, Lemma~\ref{l.harmonic}, to $W$.

\medskip\noindent\textbf{Step 1:} We prove that $v$ extends to a harmonic function $W\colon\Z^{d+1}\to\R$. Let $\mathcal T_n\colon\R\to\R$ denote the Chebyshev polynomial of degree $n$, defined by $\mathcal T_n(\cos\theta)=\cos(n\theta)$. We evaluate it at operators on functions $\Z^d\to\R$, writing $I$ for the identity, and define $W\colon\Z^{d+1}\to\R$ by
\begin{equation*}
W(x,n)\coloneqq\bigl(\mathcal T_{|n|}(I-d\Delta)v\bigr)(x)\qquad(x\in\Z^d,\ n\in\Z)\,.
\end{equation*}
Then $W(\cdot,0)=v$. The recurrence $\mathcal T_{n+1}(t)+\mathcal T_{n-1}(t)=2t\,\mathcal T_n(t)$, together with $\mathcal T_1(t)=t$ at $n=0$, yields
\begin{equation}\label{e.chebyshev}
W(\cdot,n+1)+W(\cdot,n-1)-2W(\cdot,n)=-2d\Delta W(\cdot,n)\,,
\end{equation}
which says that $W$ is harmonic on $\Z^{d+1}$.

\medskip\noindent\textbf{Step 2:} We prove that
\begin{equation}\label{e.extension}
\|W\|_{\underline L^2(Q_N^{d+1})}\leq K^r\|v\|_{L^\infty(Q_{3N}^d)}\,,\qquad\text{where }K\coloneqq3^{\nf d2}\cdot4(d+1)\,.
\end{equation}
Fix $n$ with $|n|\leq N$, and define $\varphi\colon\R\to\R$ by
\begin{equation*}
\varphi(s)\coloneqq\|\mathcal T_{|n|}(I+s\Delta)v\|_{\underline L^2(Q_N^d)}^2\,.
\end{equation*}
Since $\Delta^rv$ vanishes on $Q_{2N}^d$ and $\mathcal T_{|n|}$ has degree $|n|\leq N$, the values of the function $\mathcal T_{|n|}(I+s\Delta)v\colon\Z^d\to\R$ on $Q_N^d$ are polynomials in $s$ of degree less than $r$. So $\varphi$ is a polynomial of degree at most $2r-2$. Since $\mathcal T_{|n|}(I+s\Delta)$ is a contraction of $\ell^2(\Z^d)$ for $0\leq s\leq1$, and the values of $\mathcal T_{|n|}(I+s\Delta)v$ on $Q_N^d$ depend only on $v$ on $Q_{2N}^d$, we have
\begin{equation*}
\varphi(s)\leq\frac{\|v\|_{L^2(Q_{2N}^d)}^2}{|Q_N^d|}\leq3^d\|v\|_{L^\infty(Q_{3N}^d)}^2\qquad(0\leq s\leq1)\,.
\end{equation*}
The Remez inequality, Lemma~\ref{l.remezBLMS}, with $J=[-d,1]$ and $E=[0,1]$, extends this bound to $s=-d$ at the cost of a factor $(4(d+1))^{2r-2}$:
\begin{equation*}
\|W(\cdot,n)\|_{\underline L^2(Q_N^d)}^2=\varphi(-d)\leq\bigl(4(d+1)\bigr)^{2r-2}3^d\|v\|_{L^\infty(Q_{3N}^d)}^2\leq K^{2r}\|v\|_{L^\infty(Q_{3N}^d)}^2\,.
\end{equation*}
Averaging over the slices $|n|\leq N$ gives~\eqref{e.extension}.

\medskip\noindent\textbf{Step 3:} We conclude. Large-scale analyticity on $\Z^{d+1}$, Lemma~\ref{l.harmonic}, applied to $W$ gives a real polynomial $\tilde P\colon\R^{d+1}\to\R$ of total degree at most $m$ such that
\begin{equation*}
\|W-\tilde P\|_{L^\infty(Q_{\rho N}^{d+1})}\leq\bigl(2K^2\bigr)^{-m}\|W\|_{\underline L^2(Q_N^{d+1})}\leq2^{-m}\|v\|_{L^\infty(Q_{3N}^d)}\,,
\end{equation*}
by~\eqref{e.extension} and $r\leq2m$. Then $P\coloneqq\tilde P(\cdot,0)$ is the desired polynomial.
\end{proof}

The bound~\eqref{e.extension} is a discrete analogue of the estimate of~\citet{Kounchev} and~\citet[Theorem~1]{Pavlovic} for the iterated Laplacians of a polyharmonic function. Our proof of~\eqref{e.extension} uses Chebyshev polynomials of a transition operator, as~\citet{Carne} did in his proof of a Gaussian upper bound on the transition probabilities of reversible random walks.

\section{Proof of Theorem~\ref{t.growth}}\label{s.proof}

We use the discrete Remez inequality for polynomials of several variables; see, for example, the proof of~\citet[Proposition~A.3]{BCG}.

\begin{lemma}\label{l.remez}
Let $1\leq R\leq N$ and $0<\delta\leq1$, let $A\subset Q_R$ contain at least a fraction $\delta$ of the sites of $Q_R$, and let $F\colon\R^d\to\R$ be a real polynomial of degree at most $k$ in each variable, where $k\leq2^{-d-1}\delta R$. Then the maximum of $|F|$ over the real cube $[-3N,3N]^d$ satisfies
\begin{equation*}
\max_{[-3N,3N]^d}|F|\leq\left(\frac{C(d)N}{\delta R}\right)^{dk}\max_A|F|\,.
\end{equation*}
\end{lemma}

Let $\rho\leq1$ be given by Lemma~\ref{l.analyticity}, and let $m\coloneqq\lfloor p\rho N/2\rfloor$ and $k\coloneqq\lfloor\gamma m\rfloor$ for a small constant $\gamma=\gamma(d,p)>0$ selected below. Let $b\coloneqq(\rho\log2)/8$. We take $N_0$ at least the threshold given by Lemma~\ref{l.analyticity}, and so large that $m\geq p\rho N/4\geq1$ for $N\geq N_0$. Then $e^{bpN}\leq2^{m/2}$ and $k+1>\gamma p\rho N/4$, so it suffices to show that the set
\begin{equation*}
E\coloneqq\{x\in Q_{3N}:|u(x)|>e^{bpN}\}
\end{equation*}
has at least $(k+1)^d$ sites. Suppose, for contradiction, that $|E|<(k+1)^d$.

\medskip\noindent\textbf{Step 1:} We prove that there is a polynomial $F$ of degree at most $k$ in each variable, with $F=0$ on $E$ and $\max_{[-3N,3N]^d}|F|=1$, such that $Fu$ is close to a polynomial $P$ of total degree at most $m$,
\begin{equation}\label{e.FuP}
\|Fu-P\|_{L^\infty(Q_{\rho N})}\leq2^{-m/2}\,,
\end{equation}
and $|F|$ is small at few sites of $Q_{\rho N}$,
\begin{equation}\label{e.fewsmall}
\bigl|\{x\in Q_{\rho N}:|F(x)|\leq2^{-m/2}\}\bigr|<\frac{p|Q_{\rho N}|}{16}\,.
\end{equation}
The real polynomials on $\R^d$ of degree at most $k$ in each variable form a space of dimension $(k+1)^d>|E|$, so one of them, $F$, vanishes on $E$ and satisfies $\max_{[-3N,3N]^d}|F|=1$. Then $|Fu|\leq e^{bpN}\leq2^{m/2}$ on $Q_{3N}$. Since $u$ is harmonic and $F$ has total degree at most $dk$, an induction on the degree of $F$ shows that $\Delta^{dk+1}(Fu)=0$. Selecting $\gamma$ small enough gives $dk\leq m$, so Lemma~\ref{l.analyticity} with $r=dk+1$ gives~\eqref{e.FuP}. If~\eqref{e.fewsmall} failed, then Lemma~\ref{l.remez} with $R=\rho N$, $\delta=p/16$ and $A$ the set in~\eqref{e.fewsmall} would give
\begin{equation*}
1=\max_{[-3N,3N]^d}|F|\leq\left(\frac{16C(d)}{p\rho}\right)^{dk}2^{-m/2}<1\,,
\end{equation*}
a contradiction. Here, selecting $\gamma$ small enough makes the lemma applicable and gives the last inequality in the above display.

\medskip\noindent\textbf{Step 2:} We conclude. The restriction of $FP$ to a line is a polynomial in one variable of degree at most $m+k$, and selecting $\gamma$ small enough gives $m+k<3p\rho N/4$. So it suffices to find a line on which $FP$ changes sign across more than $3p\rho N/4$ edges. If $x\in Q_{\rho N}$ satisfies $|F(x)|>2^{-m/2}$ and $|u(x)|\geq1$, then, by~\eqref{e.FuP},
\begin{equation*}
|F(x)u(x)|>2^{-m/2}\geq|F(x)u(x)-P(x)|\,,
\end{equation*}
so $P(x)$ has the sign of $F(x)u(x)$, and $F(x)P(x)$ has the sign of $u(x)$. Thus $FP$ changes sign across every edge of $Q_{\rho N}$ that satisfies~\eqref{e.pair} and has no endpoint where $|F|\leq2^{-m/2}$. Hence, by assumption and~\eqref{e.fewsmall}, $FP$ changes sign across more than a fraction $3p/4$ of the edges of $Q_{\rho N}$ in some coordinate direction, and so across more than $3p\rho N/4$ edges of one line in this direction.\qed

\section{Periodic networks}\label{s.periodic}

In this section, we prove Theorem~\ref{t.periodic}. Recall its setting: $G$ is a connected, locally finite graph with vertex set $V$ and conductances $\bfa$, and $\Z^d$ acts freely on $V$ by bijections $\tau_x$ that preserve the conductances, with finitely many orbits. On each orbit, a harmonic function on $G$ does not in general solve the Laplace equation of $\Z^d$, but it solves a difference equation with constant coefficients, given by the determinant $S$ of a matrix of difference operators. This equation takes the place of the Laplace equation in the proof of Theorem~\ref{t.growth}. The main new ingredient is large-scale analyticity for $S$, Lemma~\ref{l.analyticityS}. To prove it, we show in Lemma~\ref{l.scalar} that $S$ is self-adjoint and coercive. This lets us turn solutions of $S^rv=0$ into polyharmonic functions in more dimensions, to which the proof of Lemma~\ref{l.analyticity} applies.

We first derive the equation on each orbit. Let $\ell$ be the number of orbits, and choose a vertex $o_i$ in the $i$th orbit for $i=1,\dots,\ell$. Since the action is free, every vertex is $(x,i)\coloneqq\tau_xo_i$ for exactly one $x\in\Z^d$ and one $i$. For $u\colon V\to\R$, define $u_i\colon\Z^d\to\R$ by $u_i(x)\coloneqq u((x,i))$. Since $\tau_x$ preserves the conductances, $\bfa((x,i),(x+y,i'))=\bfa((0,i),(y,i'))$. So $u$ is harmonic if and only if $\sum_{i'}H_{ii'}u_{i'}=0$ for every $i$, where $H$ is the $\ell\times\ell$ matrix of operators on functions $f\colon\Z^d\to\R$ with entries
\begin{equation*}
H_{ii'}f(x)\coloneqq \bfa((0,i))\delta_{ii'}f(x)-\sum_y \bfa((0,i),(y,i'))f(x+y)\,.
\end{equation*}
The entries of $H$ are finite linear combinations of translations, so they commute, and $S\coloneqq\det H$ is again such a combination. Let $\adj H$ be the adjugate of $H$, the transpose of its matrix of cofactors. Since the entries of $H$ commute, $(\adj H)H=SI$, as for matrices over a commutative ring. Applying this to $(u_1,\dots,u_\ell)$ shows that $Su_i=0$ for every $i$ if $u$ is harmonic. We also need that $S$ is self-adjoint and coercive.

\begin{lemma}\label{l.scalar}
The operator $S$ is self-adjoint on $\ell^2(\Z^d)$ and coercive: for some $\theta>0$,
\begin{equation*}
\langle Sg,g\rangle_{\ell^2(\Z^d)}\geq \theta\sum_{j=1}^d\|g(\cdot+e_j)-g\|_{\ell^2(\Z^d)}^2\qquad\text{for every finitely supported }g\colon\Z^d\to\R\,.
\end{equation*}
\end{lemma}

\begin{proof}
We pass to Fourier series, which turn $S$ into multiplication by the determinant of a Hermitian matrix, and bound this determinant from below using paths in $G$. For $z\in\T^d$, where $\T\coloneqq\{\zeta\in\C:|\zeta|=1\}$, let $H(z)$ be the $\ell\times\ell$ matrix obtained from $H$ by replacing each translation $f\mapsto f(\cdot+y)$ by the number $z^y\coloneqq z_1^{y_1}\cdots z_d^{y_d}$.

\medskip\noindent\textbf{Step 1:} We prove that $H(z)$ is Hermitian and that there is $c>0$, independent of $z$, such that
\begin{equation*}
v^*H(z)v\geq c\biggl(\sum_{j=1}^d|z_j-1|^2|v_1|^2+\sum_{i=2}^\ell|v_i-v_1|^2\biggr)\qquad\text{for every }v\in\C^\ell\,.
\end{equation*}
Since the conductances are symmetric and $\tau_{-y}$ preserves them, $\bfa((0,i'),(y,i))=\bfa((0,i),(-y,i'))$. Since also $\overline{z^y}=z^{-y}$, the matrix $H(z)$ is Hermitian, and expanding the squares shows that
\begin{equation*}
v^*H(z)v=\frac12\sum_{i,i',y}\bfa((0,i),(y,i'))\,|v_i-z^yv_{i'}|^2\,.
\end{equation*}
So, along an edge between $(x,i)$ and $(x+y,i')$, the increment of the function $(x,i)\mapsto z^xv_i$ has modulus $|v_i-z^yv_{i'}|$, which is at most $(2v^*H(z)v)^{1/2}$ divided by the square root of the total conductance between the endpoints. Since $G$ is connected, we can fix paths in $G$ from $(0,1)$ to each $(0,i)$ and to each $(e_j,1)$, along which the total increments of this function are $v_i-v_1$ and $(z_j-1)v_1$. Adding the bounds along each path, squaring, and summing over the paths gives the claim, with a constant $c$ that does not depend on $z$ because the paths do not.

\medskip\noindent\textbf{Step 2:} We prove that $\det H(z)\geq c^\ell\sum_j|z_j-1|^2$. By Step~1, $H(z)$ is positive semidefinite, so $\det H(z)\geq0$, which proves the claim if $z=(1,\dots,1)$. Otherwise, $z_j\neq1$ for some $j$, and the lower bound of Step~1 is $c|w|^2$ in the coordinates $w=\bigl((\sum_j|z_j-1|^2)^{1/2}v_1,v_2-v_1,\dots,v_\ell-v_1\bigr)$, whose change of variables has determinant $(\sum_j|z_j-1|^2)^{1/2}$. So the matrix of the quadratic form in these coordinates has all eigenvalues at least $c$, hence determinant at least $c^\ell$, and $\det H(z)$ is $\sum_j|z_j-1|^2$ times this determinant.

\medskip\noindent\textbf{Step 3:} We conclude. For finitely supported $g\colon\Z^d\to\R$, the Fourier series $\hat g(z)\coloneqq\sum_xg(x)z^{-x}$ turns $g(\cdot+y)$ into $z^y\hat g$. Since composing translations adds their shifts, it turns $Sg$ into $\det H(z)\,\hat g$. By Parseval's identity, with integrals over $\T^d$ against the uniform probability measure, it follows that $\langle Sf,g\rangle=\int\det H(z)\,\hat f\,\overline{\hat g}=\langle f,Sg\rangle$ for finitely supported $f$ and $g$, since $\det H(z)$ is real, and, by Step~2,
\begin{equation*}
\langle Sg,g\rangle=\int\det H(z)\,|\hat g|^2\geq c^\ell\int\sum_{j=1}^d|z_j-1|^2|\hat g|^2=c^\ell\sum_{j=1}^d\|g(\cdot+e_j)-g\|^2\,.
\end{equation*}
Since $S$ is a finite linear combination of translations, it is bounded on $\ell^2(\Z^d)$, and finitely supported functions are dense in $\ell^2(\Z^d)$. So $S$ is self-adjoint, and the lemma holds with $\theta=c^\ell$.
\end{proof}

The proof of Theorem~\ref{t.growth} relies on Lemma~\ref{l.analyticity}, large-scale analyticity for solutions of $\Delta^rv=0$. Using Lemma~\ref{l.scalar}, we prove the same estimate for solutions of $S^rv=0$.

\begin{lemma}\label{l.analyticityS}
Lemma~\ref{l.analyticity} holds with $S$ in place of $\Delta$, with $\rho$ and $N_0$ depending on $S$.
\end{lemma}

\begin{proof}
We turn solutions of $S^rv=0$ into polyharmonic functions on $\Z^{d+q}$, for a suitable $q$, and apply the proof of Lemma~\ref{l.analyticity} there. We write $\Delta_d$ and $\Delta_{d+q}$ for the Laplacians on $\Z^d$ and $\Z^{d+q}$, defined as in the introduction. Replacing $S$ by $S/(2\theta)$ does not change the solutions of $S^rv=0$, so we may assume that $\theta=1/2$ in Lemma~\ref{l.scalar}. Since $\langle-d\Delta_d g,g\rangle=\frac12\sum_j\|g(\cdot+e_j)-g\|^2$, Lemma~\ref{l.scalar} gives $\langle(S+d\Delta_d)g,g\rangle\geq0$ for finitely supported $g$, and hence, by density, $S+d\Delta_d\geq0$ on $\ell^2(\Z^d)$. Fix a positive integer $q\geq\|S+d\Delta_d\|_{\ell^2(\Z^d)}$, so that $0\leq I-(S+d\Delta_d)/q\leq I$. As in the proof of Lemma~\ref{l.analyticity}, let $\mathcal T_n$ be the Chebyshev polynomial of degree $n$, given by $\mathcal T_n(\cos t)=\cos(nt)$, and define $\tilde v\colon\Z^{d+q}\to\R$ by
\begin{equation*}
\tilde v(x,z)\coloneqq\bigl(\mathcal T_{|z_1+\cdots+z_q|}\bigl(I-(S+d\Delta_d)/q\bigr)v\bigr)(x)\qquad(x\in\Z^d,\ z\in\Z^q)\,.
\end{equation*}
Then $\tilde v(\cdot,0)=v$. Below, $S$ and $\Delta_d$ act on functions on $\Z^{d+q}$ through the coordinates $x_1,\dots,x_d$.

\medskip\noindent\textbf{Step 1:} We prove that
\begin{equation*}
\bigl((d+q)\Delta_{d+q}\bigr)^r\tilde v=(-S)^r\tilde v\,.
\end{equation*}
As in~\eqref{e.chebyshev}, the recurrence of the Chebyshev polynomials shows that the second difference of $\tilde v$ in each of the coordinates $z_1,\dots,z_q$ is $-(2/q)(S+d\Delta_d)\tilde v$. The second differences in the coordinates $x_1,\dots,x_d$ add up to $2d\Delta_d\tilde v$. So $(d+q)\Delta_{d+q}\tilde v=-S\tilde v$. Iterating, since $S$ commutes with $\Delta_{d+q}$, gives the claim.

\medskip\noindent\textbf{Step 2:} We prove that $\tilde v$ is polyharmonic in a box smaller than $Q_{3N}$ by a constant factor, with controlled average size there. More precisely, there are an integer $M\geq1$ and $C\geq1$, depending only on $S$, such that, with $N'\coloneqq\lfloor N/M\rfloor$,
\begin{equation*}
\Delta_{d+q}^r\tilde v=0\quad\text{in }Q_{3N'}^{d+q}\qquad\text{and}\qquad\|\tilde v\|_{\underline L^2(Q_{3N'}^{d+q})}\leq C\|v\|_{L^\infty(Q_{3N})}\,.
\end{equation*}
The operator $\mathcal T_n(I-(S+d\Delta_d)/q)$ is a polynomial of degree $n$ in an operator of finite range, and it commutes with $S$. So $\tilde v(x,z)$ and $S^r\tilde v(x,z)$ depend only on the values of $v$ and of $S^rv$ within distance proportional to $|z_1+\cdots+z_q|$ of $x$. For $M$ large, these values lie in $Q_{3N}$ when $(x,z)\in Q_{3N'}^{d+q}$, and $S^rv=0$ there, so Step~1 gives the first claim. Since $0\leq I-(S+d\Delta_d)/q\leq I$, each $\mathcal T_n(I-(S+d\Delta_d)/q)$ is a contraction of $\ell^2(\Z^d)$, and averaging over $z$ gives the second claim.

\medskip\noindent\textbf{Step 3:} We conclude by applying the proof of Lemma~\ref{l.analyticity} to $\tilde v$ on $\Z^{d+q}$ at scale $N'$. Its Step~2 uses the size of $\tilde v$ only through $\|\tilde v\|_{\underline L^2(Q_{3N'}^{d+q})}$, which our Step~2 bounds by $C\|v\|_{L^\infty(Q_{3N})}$, and applying Lemma~\ref{l.harmonic} with a smaller $\alpha$ in Step~3 of that proof absorbs the factor $C$. Restricting the resulting polynomial to $z=0$ gives the lemma.
\end{proof}

With Lemma~\ref{l.analyticityS}, we can now prove Theorem~\ref{t.periodic}.

\begin{proof}[Proof of Theorem~\ref{t.periodic}]
Choose $o$ as the representative of its orbit, so that $u(\tau_xo)=u_i(x)$ for the corresponding $i$. Since $u$ is harmonic, $Su_i=0$. We repeat the proof of Theorem~\ref{t.growth} for $u_i$, with $S$ in place of $\Delta$ and Lemma~\ref{l.analyticityS} in place of Lemma~\ref{l.analyticity}. The induction in Step~1 that gives $\Delta^{dk+1}(Fu)=0$ uses only that $\Delta$ is a finite linear combination of translations, so it also gives $S^{dk+1}(Fu_i)=0$.
\end{proof}

\section{Uniformly elliptic conductances}\label{s.random}

In this section, we prove Theorems~\ref{t.random} and~\ref{t.counter}. In Section~\ref{s.planar}, we use planarity to show that every stationary $\bfa$-harmonic function on $\Z^2$ with $\lambda^{-1}\leq\bfa\leq\lambda$ is constant, and in Section~\ref{s.examples}, we construct nonconstant examples with $\bfa\leq1$ and with $\bfa\geq1$. In Section~\ref{s.higher}, for $d\geq3$, we glue together harmonic functions supported on planes into a nonconstant stationary function that is nearly harmonic, and then choose conductances $\bfa$ close to $1$ for which it is $\bfa$-harmonic.

\subsection{Constancy on the square lattice}\label{s.planar}

We prove Theorem~\ref{t.random}(a) by contradiction. If $h$ is not constant, then $t_n\coloneqq\max_{Q_{2n}}|h|$ tends to infinity. We show that $|h|>t_n/2$ on a positive fraction of $Q_{4n}$ for all large $n$, which contradicts stationarity. This uses two ingredients. The first is topological. By the maximum principle, the sets $\{h>t\}$ and $\{h<-t\}$ have no bounded components. Lemma~\ref{l.planar} combines this with planarity and stationarity to show that the set $\{|h|\leq t\}$ between them has no saddles. That is, no finite region has $\{h>t\}$ and $\{h<-t\}$ alternating four times around it. By the maximum principle, $|h|>t_n$ at some site just outside $Q_{2n}$. Near this site, the absence of saddles gives a region $U$ of size comparable to $n$ on which $|h|\leq t_n$. It borders only one of $\{h>t_n\}$ and $\{h<-t_n\}$, say $\{h>t_n\}$. The second ingredient is analytic: since $h\geq-t_n$ on $U$, Lemma~\ref{l.hitting}, elliptic regularity at the boundary, shows that $h>t_n/2$ on a square of size comparable to $n$ in $U$ near $\{h>t_n\}$.

Throughout, $h$ and $\bfa$ are as in Theorem~\ref{t.random}(a), and we write $Q_s(x)\coloneqq x+Q_s$.

\begin{lemma}\label{l.planar}
Almost surely, the following holds for every $t>0$. Let $B\subseteq\{|h|\leq t\}$ be finite and connected, and let $U_1,U_2,U_3\subseteq\{|h|\leq t\}\setminus B$ be disjoint, finite and connected, each adjacent to $B$. Then some $U_i$ is not adjacent to both $\{h>t\}$ and $\{h<-t\}$.
\end{lemma}

\begin{figure}[ht]
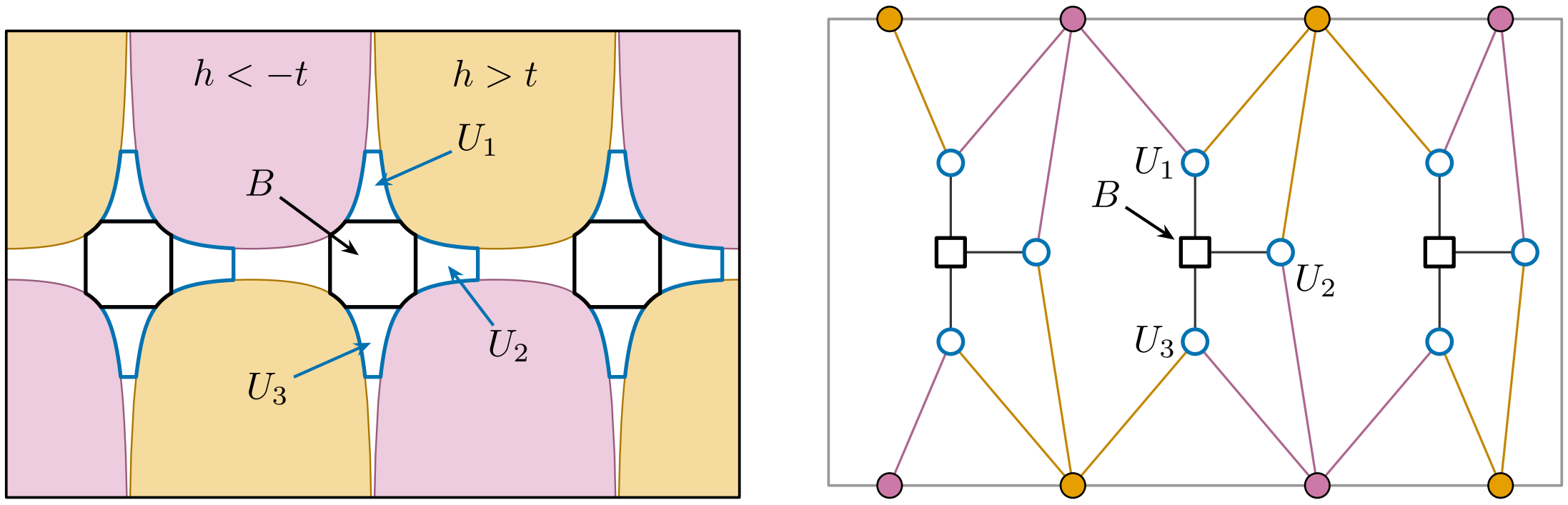

\centering
\definecolor{confU}{HTML}{56B4E9}
\definecolor{confP}{HTML}{E69F00}
\definecolor{confM}{HTML}{CC79A7}
\definecolor{confArc}{HTML}{0B3C5D}
\definecolor{confB}{HTML}{000000}%
\definecolor{confArm}{HTML}{0072B2}%

\caption{Left: the sets $\{h>t\}$ (orange), $\{h<-t\}$ (purple) and $\{|h|\leq t\}$ (white) for the harmonic function $h(x,y)=\sin x\sinh y$ and small $t>0$, with sets $B$ (outlined in black) and $U_1,U_2,U_3$ (outlined in blue) at each saddle point. Right: the contracted graph. Each $B$ is a black square, each $U_i$ a blue circle, and each component of $\{h>t\}$ or of $\{h<-t\}$ an orange or purple disk on the frame; each $U_i$ is joined to its $B$ and to one component of each sign.}
\label{f.contract}
\end{figure}

\begin{proof}
Figure~\ref{f.contract} (left) shows what the lemma rules out. The function $\sin x\sinh y$ is harmonic but not stationary. At each of its saddle points, there are sets $B,U_1,U_2,U_3$ as in the lemma with every $U_i$ adjacent to both $\{h>t\}$ and $\{h<-t\}$. Step~1 shows that a square of side $n$ contains at most of order $n$ disjoint families of such sets. This step is deterministic: it uses only the maximum principle and planarity. Step~2 shows that a stationary function has no such families at all.

\medskip\noindent\textbf{Step 1:} We prove that, for every $t>0$ and $n\geq1$, the square $Q_n$ contains at most $16n$ disjoint such families. Here a set contains a family if it contains the sets $B$ and $U_i$ and, for each $i$, a site of $\{h>t\}$ and a site of $\{h<-t\}$ adjacent to $U_i$. Let $m$ be the number of disjoint such families in $Q_n$. By the maximum principle, every component of $\{h>t\}$ and of $\{h<-t\}$ is infinite. So every component of their intersections with $Q_n$ contains one of the $8n$ sites of $Q_n$ that have a neighbor outside $Q_n$, just as each orange or purple region in Figure~\ref{f.contract} (left) reaches the frame. Contract each of these components, and each set $B$ and $U_i$, to a vertex. Join each $U_i$ to its $B$ and to one adjacent component of each sign. This gives a planar bipartite graph, as on the right of Figure~\ref{f.contract}, with $9m$ edges and at most $4m+8n$ vertices. Euler's formula for planar graphs gives $9m\leq2(4m+8n)$, that is, $m\leq16n$.

\medskip\noindent\textbf{Step 2:} We prove the lemma by comparing Step~1 with stationarity. If such families occurred with positive probability, then a square of side $n$ would contain on average of order $n^2$ disjoint ones, more than Step~1 allows. Precisely, if the lemma fails, then there are $t>0$ and $\ell\geq1$ such that $Q_\ell$ contains such a family with some probability $p>0$. For $n\geq\ell$, the square $Q_n$ contains at least $(n/(3\ell))^2$ disjoint translates of $Q_\ell$. By stationarity, the expected number of disjoint such families in $Q_n$ is therefore at least $p(n/(3\ell))^2$, which exceeds $16n$ for large $n$.
\end{proof}

Uniform ellipticity, which is necessary by Theorem~\ref{t.random}(b), enters through the following standard estimate, a form of Hölder continuity at the boundary. Here, a function $w$ is \textbf{$\bfa$-subharmonic} on a set $D$ if $\sum_{|y-z|=1}\bfa(z,y)(w(y)-w(z))\geq0$ for every $z\in D$.

\begin{lemma}[{see \citealp[Lemma~2.11]{Chelkak}}]\label{l.hitting}
There are $C,\beta>0$, depending only on $\lambda$, with the following property. Let $x\in\Z^2$, let $1\leq s\leq R$, and let $S\subseteq\Z^2$ be infinite and connected with $S\cap Q_s(x)\neq\emptyset$. If $w\colon\Z^2\to[0,1]$ vanishes on $S$ and is $\bfa$-subharmonic on $Q_R(x)\setminus S$, then $w\leq C(s/R)^\beta$ on $Q_s(x)$.
\end{lemma}

With Lemmas~\ref{l.planar} and~\ref{l.hitting}, we can now prove Theorem~\ref{t.random}(a).

\begin{figure}[ht]
\centering
\definecolor{annArc}{HTML}{0B3C5D}%
\definecolor{annU}{HTML}{56B4E9}%
\definecolor{annAp}{HTML}{E69F00}%
\definecolor{annAm}{HTML}{CC79A7}%
\begin{tikzpicture}[x=1.6cm,y=1.6cm,line join=round,
  lab/.style={font=\small,inner sep=1.5pt},
  arc/.style={annArc,line width=2pt,line cap=butt,line join=miter}]
  \begin{scope}
  \clip (-1.15,-2.25) rectangle (3.75,2.25);
  \def\sqO{(0,-2.06) rectangle (4.12,2.06)}
  \def\sqM{(0.3333,-1.7267) rectangle (3.7867,1.7267)}
  \def\sqI{(0.6667,-1.3933) rectangle (3.4533,1.3933)}
  \def\sqK{(1,-1.06) rectangle (3.12,1.06)}
  \fill[black!7,even odd rule] \sqO \sqK;
  \begin{scope}
    \clip (-3,1) rectangle (5,3);
    \fill[annU!75,even odd rule] \sqI \sqK;
  \end{scope}
  \begin{scope}
    \clip (-3,-3) rectangle (5,-1);
    \fill[annU!30,even odd rule] \sqM \sqI;
  \end{scope}
  \begin{scope}
    \clip (-3,-3) rectangle (3.5,-1);
    \fill[annU!30,even odd rule] \sqO \sqM;
  \end{scope}
  \draw[black!35,line width=0.4pt] \sqO \sqM \sqI \sqK;
  \fill[white] (-1,-1) rectangle (1,1);
  \draw[black!70,line width=0.5pt,line join=miter] (-1,-1) rectangle (1,1);
  \draw[black!60,dashed,line width=0.5pt] (-2,-2) rectangle (2,2);
  \def\setAp{plot[smooth,tension=0.6] coordinates{(3.05,2.6) (2.92,2.15) (2.98,1.8) (2.8,1.5)
    (2.75,1.23) (2.58,1.0) (2.47,0.72) (2.27,0.55) (2.1,0.32) (2.03,0.05) (2.05,-0.22)
    (2.28,-0.38) (2.42,-0.66) (2.78,-0.84) (3.02,-1.13) (3.25,-1.37) (3.6,-1.44) (4.2,-1.45)
    (4.9,-1.46)}
    -- (4.9,2.6) -- cycle}
  \def\setAm{plot[smooth,tension=0.6] coordinates{(4.9,-1.51) (4.3,-1.52) (3.85,-1.51)
    (3.55,-1.6) (3.32,-1.56) (3.12,-1.66) (3.0,-1.79) (2.96,-1.9) (2.86,-2.06) (2.93,-2.32)
    (2.84,-2.6)} -- (4.9,-2.6) -- cycle}
  \fill[annAp!38] \setAp;
  \fill[annAm!40] \setAm;
  \draw[annAp!80!black,line width=0.5pt] \setAp;
  \draw[annAm!80!black,line width=0.5pt] \setAm;
  \draw[arc] (0.8333,1) -- (0.8333,1.2267) -- (2.74,1.2267);
  \draw[arc] (0.5,-1) -- (0.5,-1.56) -- (3.3,-1.56);
  \draw[arc] (0.1667,-1) -- (0.1667,-1.8933) -- (2.95,-1.8933);
  \draw[black,line width=0.5pt,dash pattern=on 2pt off 1.5pt] (2.46,1.0667) rectangle (2.78,1.3867);
  \filldraw[fill=annU,draw=black,line width=0.4pt] (2.565,1.1717) rectangle (2.675,1.2817);
  \fill (2.62,1.2267) circle (0.9pt);
  \draw[black!70,line width=0.5pt,-{Stealth[length=4pt,width=3pt]}] (2.2,1.75) -- (2.58,1.27);
  \node[lab,anchor=south] at (2.2,1.72) {$v$};
  \end{scope}
  \fill (2.06,0) circle (1.7pt);
  \node[lab,right=1pt] at (2.06,0) {$x$};
  \node[lab] at (-0.45,0.55) {$Q_n$};
  \node[lab,anchor=east] at (-1.17,2.0) {$Q_{2n}$};
  \node[lab] at (3.25,0.6) {$A_+$};
  \draw[black!70,line width=0.5pt,-{Stealth[length=4pt,width=3pt]}] (0.2,1.62) -- (0.74,1.3);
  \node[lab,anchor=east,inner sep=1pt] at (0.2,1.62) {$U$};
  \node[lab] at (3.4,-2.0) {$A_-$};
\end{tikzpicture}%
\caption{Three annuli (gray) about a site $x\in\partial Q_{2n}$ in $A_+\cup A_-$ and, in each, an arc (thick) from $Q_n$ to $A_+\cup A_-$ and the component (blue) of the annulus minus $Q_n\cup A_+\cup A_-$ that contains it. On the arc in $U$, a site $v$ with the squares $Q_{\delta n/2}(v)$ (filled) and $Q_{n/12}(v)$ (dashed), not to scale.}
\label{f.annuli}
\end{figure}
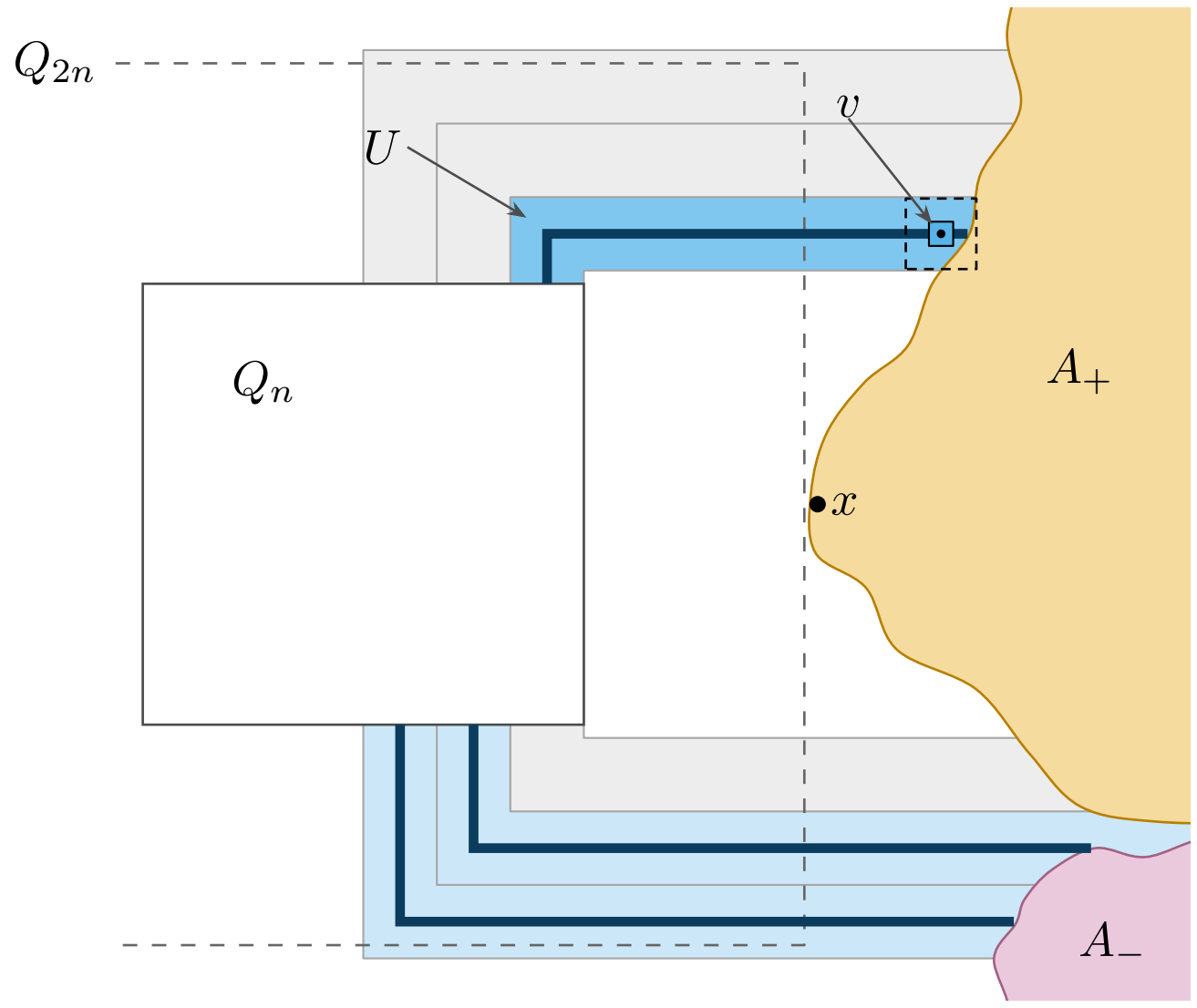

\begin{proof}[Proof of Theorem~\ref{t.random}(a)]
We may suppose that $h$ is almost surely not constant. Then $t_n\coloneqq\max_{Q_{2n}}|h|$ tends to infinity by the Liouville theorem \citep[see, for example,][Section~2.5 and Exercise~2.43]{LyonsPeres}. So, by stationarity, it suffices to show that $|h|>t_n/2$ on at least $c(\lambda)n^2$ sites of $Q_{4n}$ for all large $n$. We fix a realization of $h$ and $\bfa$ for which $h$ is not constant and the conclusion of Lemma~\ref{l.planar} holds, let $n$ be large, and let $A_\pm\coloneqq\{\pm h>t_n\}$ and $A\coloneqq A_+\cup A_-$.

\medskip\noindent\textbf{Step 1:} We find the region $U$ shown in Figure~\ref{f.annuli}. By the maximum principle, some $x\in\partial Q_{2n}$ lies in $A$, and the component of $A$ that contains $x$ is infinite. Figure~\ref{f.annuli} shows the three annuli (gray) of sites $z$ with $n+jn/3<|z-x|_\infty<n+(j+1)n/3$, for $j\in\{0,1,2\}$. Each annulus contains, in its middle third, an arc (thick) from $Q_n$ to $A$ whose interior lies in a component (blue) of the annulus minus $Q_n\cup A$. By Lemma~\ref{l.planar} with $t=t_n$ and $B=Q_n$, one of these three components, which we call $U$ (dark blue), is adjacent to only one of $A_+$ and $A_-$, say to $A_+$.

\medskip\noindent\textbf{Step 2:} We prove that $h>t_n/2$ on a square in $U$ with at least $c(\lambda)n^2$ sites. Let $\delta=\delta(\lambda)>0$ be small, and let $v$ be a site of the arc in $U$ with $\min_{z\in A}|v-z|_\infty=\lfloor\delta n\rfloor$ (Figure~\ref{f.annuli}). Then $Q_{\delta n/2}(v)\subseteq U$, and $Q_{n/12}(v)$ is disjoint from $Q_n$ and lies in the annulus that contains $U$. Let $w\coloneqq(t_n-h)_+/(2t_n)$ on $U$ and $w\coloneqq0$ elsewhere. Then $0\leq w\leq1$. Since $(t_n-h)_+$ is $\bfa$-subharmonic and vanishes on $A_+$, the only part of $A$ adjacent to $U$, the function $w$ is $\bfa$-subharmonic on $Q_{n/12}(v)$. We apply Lemma~\ref{l.hitting} with $v$ in place of $x$, $s=\delta n$ and $R=n/12$, and with $S$ a component of $A$ that meets $Q_{\delta n}(v)$. It gives $w\leq C(12\delta)^\beta$ on $Q_{\delta n}(v)$, which is less than $\nf14$ if $\delta$ is small enough. So $h>t_n/2$ on $Q_{\delta n/2}(v)\subseteq Q_{4n}$.
\end{proof}

\subsection{Counterexamples on the square lattice}\label{s.examples}

We prove Theorem~\ref{t.random}(b) with a stationary random forest on $\Z^2$ in which every site has finitely many descendants. A similar forest appears in the counterexample of~\citet[Section~3]{ZernerMerkl} for planar random walks in random environment.

\begin{proof}[Proof of Theorem~\ref{t.random}(b)]
Divide $\Z^2$ into the $2\times2$ blocks $2z+\{0,1\}^2$ with $z\in\Z^2$, and give each row and each column of blocks an independent label, uniform on $[0,1]$. Each site points to a neighbor, as in Figure~\ref{f.blocks}. In the blocks whose row label is larger than their column label (shaded), the lower sites point east and the upper sites point west. In the other blocks, the left sites point north and the right sites point south. The site that $x$ points to is its parent $p(x)$, and $\mathcal D(x)$ is the number of descendants of $x$, including $x$ itself. Let $\sigma(x)\coloneqq1$ if $x$ points east or north (filled in Figure~\ref{f.blocks}), and $\sigma(x)\coloneqq-1$ otherwise (open).

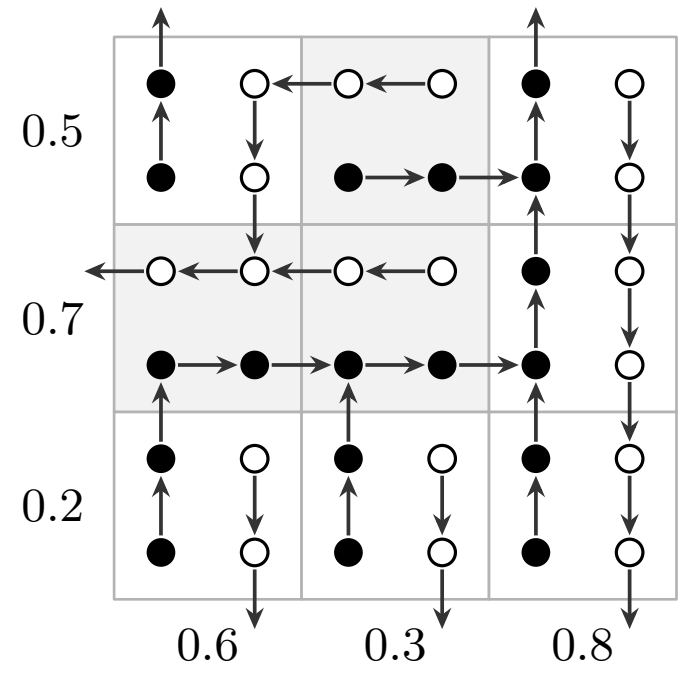
\begin{figure}[ht]
\centering
\begin{tikzpicture}[x=0.62cm,y=0.62cm,line join=round,
  lab/.style={font=\footnotesize,inner sep=1pt},
  arr/.style={black!80,line width=0.7pt,-{Stealth[length=4pt,width=3.4pt]},shorten <=3.2pt,shorten >=3.2pt},
  plus/.style={circle,fill=black,draw=black,inner sep=0pt,minimum size=5pt},
  minus/.style={circle,fill=white,draw=black,line width=0.6pt,inner sep=0pt,minimum size=5pt}]
  \draw[black!30,line width=0.5pt,fill=white] (-0.5,-0.5) rectangle (1.5,1.5);
  \draw[black!30,line width=0.5pt,fill=black!5] (-0.5,1.5) rectangle (1.5,3.5);
  \draw[black!30,line width=0.5pt,fill=white] (-0.5,3.5) rectangle (1.5,5.5);
  \draw[black!30,line width=0.5pt,fill=white] (1.5,-0.5) rectangle (3.5,1.5);
  \draw[black!30,line width=0.5pt,fill=black!5] (1.5,1.5) rectangle (3.5,3.5);
  \draw[black!30,line width=0.5pt,fill=black!5] (1.5,3.5) rectangle (3.5,5.5);
  \draw[black!30,line width=0.5pt,fill=white] (3.5,-0.5) rectangle (5.5,1.5);
  \draw[black!30,line width=0.5pt,fill=white] (3.5,1.5) rectangle (5.5,3.5);
  \draw[black!30,line width=0.5pt,fill=white] (3.5,3.5) rectangle (5.5,5.5);
  \draw[arr] (0,0) -- (0,1);
  \draw[arr] (0,1) -- (0,2);
  \draw[arr] (0,2) -- (1,2);
  \draw[arr] (0,3) -- (-1,3);
  \draw[arr] (0,4) -- (0,5);
  \draw[arr] (0,5) -- (0,6);
  \draw[arr] (1,0) -- (1,-1);
  \draw[arr] (1,1) -- (1,0);
  \draw[arr] (1,2) -- (2,2);
  \draw[arr] (1,3) -- (0,3);
  \draw[arr] (1,4) -- (1,3);
  \draw[arr] (1,5) -- (1,4);
  \draw[arr] (2,0) -- (2,1);
  \draw[arr] (2,1) -- (2,2);
  \draw[arr] (2,2) -- (3,2);
  \draw[arr] (2,3) -- (1,3);
  \draw[arr] (2,4) -- (3,4);
  \draw[arr] (2,5) -- (1,5);
  \draw[arr] (3,0) -- (3,-1);
  \draw[arr] (3,1) -- (3,0);
  \draw[arr] (3,2) -- (4,2);
  \draw[arr] (3,3) -- (2,3);
  \draw[arr] (3,4) -- (4,4);
  \draw[arr] (3,5) -- (2,5);
  \draw[arr] (4,0) -- (4,1);
  \draw[arr] (4,1) -- (4,2);
  \draw[arr] (4,2) -- (4,3);
  \draw[arr] (4,3) -- (4,4);
  \draw[arr] (4,4) -- (4,5);
  \draw[arr] (4,5) -- (4,6);
  \draw[arr] (5,0) -- (5,-1);
  \draw[arr] (5,1) -- (5,0);
  \draw[arr] (5,2) -- (5,1);
  \draw[arr] (5,3) -- (5,2);
  \draw[arr] (5,4) -- (5,3);
  \draw[arr] (5,5) -- (5,4);
  \node[plus] at (0,0) {};
  \node[plus] at (0,1) {};
  \node[plus] at (0,2) {};
  \node[minus] at (0,3) {};
  \node[plus] at (0,4) {};
  \node[plus] at (0,5) {};
  \node[minus] at (1,0) {};
  \node[minus] at (1,1) {};
  \node[plus] at (1,2) {};
  \node[minus] at (1,3) {};
  \node[minus] at (1,4) {};
  \node[minus] at (1,5) {};
  \node[plus] at (2,0) {};
  \node[plus] at (2,1) {};
  \node[plus] at (2,2) {};
  \node[minus] at (2,3) {};
  \node[plus] at (2,4) {};
  \node[minus] at (2,5) {};
  \node[minus] at (3,0) {};
  \node[minus] at (3,1) {};
  \node[plus] at (3,2) {};
  \node[minus] at (3,3) {};
  \node[plus] at (3,4) {};
  \node[minus] at (3,5) {};
  \node[plus] at (4,0) {};
  \node[plus] at (4,1) {};
  \node[plus] at (4,2) {};
  \node[plus] at (4,3) {};
  \node[plus] at (4,4) {};
  \node[plus] at (4,5) {};
  \node[minus] at (5,0) {};
  \node[minus] at (5,1) {};
  \node[minus] at (5,2) {};
  \node[minus] at (5,3) {};
  \node[minus] at (5,4) {};
  \node[minus] at (5,5) {};
  \node[lab,anchor=east] at (-0.75,0.5) {$0.2$};
  \node[lab,anchor=east] at (-0.75,2.5) {$0.7$};
  \node[lab,anchor=east] at (-0.75,4.5) {$0.5$};
  \node[lab,anchor=north] at (0.5,-0.75) {$0.6$};
  \node[lab,anchor=north] at (2.5,-0.75) {$0.3$};
  \node[lab,anchor=north] at (4.5,-0.75) {$0.8$};
\end{tikzpicture}
\caption{A window of $3\times3$ blocks, with the labels of their rows (left) and columns (bottom). In the shaded blocks, the row label is larger. Each site points to its parent; filled sites have $\sigma=1$ and open sites $\sigma=-1$.}
\label{f.blocks}
\end{figure}

\medskip\noindent\textbf{Step 1:} We prove that, almost surely, $\mathcal D(x)<\infty$ for every site $x$. As in Figure~\ref{f.blocks}, an arrow that leaves a shaded block enters the neighboring block in the same row of blocks, and an arrow that leaves an unshaded block enters the neighboring block in the same column. So the larger of the two labels of a block does not decrease along arrows. Hence the descendants of a site $x$ lie in blocks whose labels are at most the larger label of the block of $x$, which we call $m$. Almost surely, some rows of blocks above and below $x$, and some columns to its left and right, have labels larger than $m$, and these enclose all the descendants of $x$.

\medskip\noindent\textbf{Step 2:} We prove that $h\coloneqq\sigma\mathcal D^2$ is $\bfa$-harmonic for some conductances $0<\bfa\leq1$. Checking neighboring blocks, as in Figure~\ref{f.blocks}, shows that two neighbors have the same sign exactly when one of them is the parent of the other. Hence the arrows along a chain $x,p(x),p(p(x)),\dots$ all point east or north, or all point west or south, so there are no cycles. We choose the conductances
\begin{equation*}
\bfa(x,p(x))\coloneqq\frac{2\mathcal D(x)+1}{\mathcal D(p(x))^2-\mathcal D(x)^2}\qquad\text{and}\qquad\bfa(x,y)\coloneqq\frac1{\mathcal D(x)^2+\mathcal D(y)^2}\quad\text{if }\sigma(x)\neq\sigma(y)\,,
\end{equation*}
so that the current $J(x,y)\coloneqq\bfa(x,y)(h(y)-h(x))$ from $x$ to a neighbor $y$ is
\begin{equation*}
\sigma(x)\bigl(2\mathcal D(x)+1\bigr)\ \text{if }y=p(x),\qquad-\sigma(x)\bigl(2\mathcal D(y)+1\bigr)\ \text{if }p(y)=x,\qquad-\sigma(x)\ \text{if }\sigma(y)\neq\sigma(x)\,.
\end{equation*}
Since $\mathcal D(x)=1+\sum_{p(y)=x}\mathcal D(y)$, we have $\mathcal D(p(x))\geq\mathcal D(x)+1$, so $0<\bfa\leq1$, and, if $x$ has $k$ children,
\begin{equation*}
\sigma(x)\sum_{|y-x|=1}J(x,y)=2\mathcal D(x)+1-\sum_{p(y)=x}\bigl(2\mathcal D(y)+1\bigr)-(3-k)=0\,.
\end{equation*}

\medskip\noindent\textbf{Step 3:} We prove that, after a random translation, $(\bfa,h)$ is jointly stationary and $h$ is not constant. Translations by $2\Z^2$ preserve the law of the labels, so translating $\bfa$ and $h$ by a uniform element of $\{0,1\}^2$, independent of the labels, makes $(\bfa,h)$ jointly stationary. Since every block contains sites of both signs, $h$ is not constant.

\medskip\noindent\textbf{Step 4:} We construct the example with conductances at least $1$. By Step~2, $|J|\geq1$. Let $\varphi(r)\coloneqq r/(1+2|r|)$, which is strictly increasing with values in $(-\frac12,\frac12)$, let $\tilde h\coloneqq\varphi\circ h$, and let $\tilde\bfa(x,y)\coloneqq J(x,y)/(\tilde h(y)-\tilde h(x))$. Since $\varphi$ is increasing and $|\tilde h(y)-\tilde h(x)|<1$, we have $\tilde\bfa>|J|\geq1$. Moreover, $\sum_{|y-x|=1}\tilde\bfa(x,y)(\tilde h(y)-\tilde h(x))=\sum_{|y-x|=1}J(x,y)=0$, so $\tilde h$ is $\tilde\bfa$-harmonic. The pair $(\tilde\bfa,\tilde h)$ is jointly stationary, and $\tilde h$ is bounded and not constant.
\end{proof}

\subsection{Counterexamples in dimensions three and higher}\label{s.higher}

To prove Theorem~\ref{t.counter}, we glue together pieces of the harmonic function $(3+2\sqrt2)^{x_3}(-1)^{x_1}\mathbf 1\{x_1=x_2\}$ and of its images under the symmetries of $\Z^3$. The resulting function $h$ is nearly harmonic: every site $x$ has a neighbor $p(x)$, its parent, with $|h(p(x))|\geq(3+2\sqrt2)|h(x)|$, and $|\Delta h(x)|$ is small compared with $|h(p(x))|$. Changing the conductances on the edges $\{x,p(x)\}$ then gives conductances $\bfa$ for which $h$ is $\bfa$-harmonic.

\begin{figure}[t]
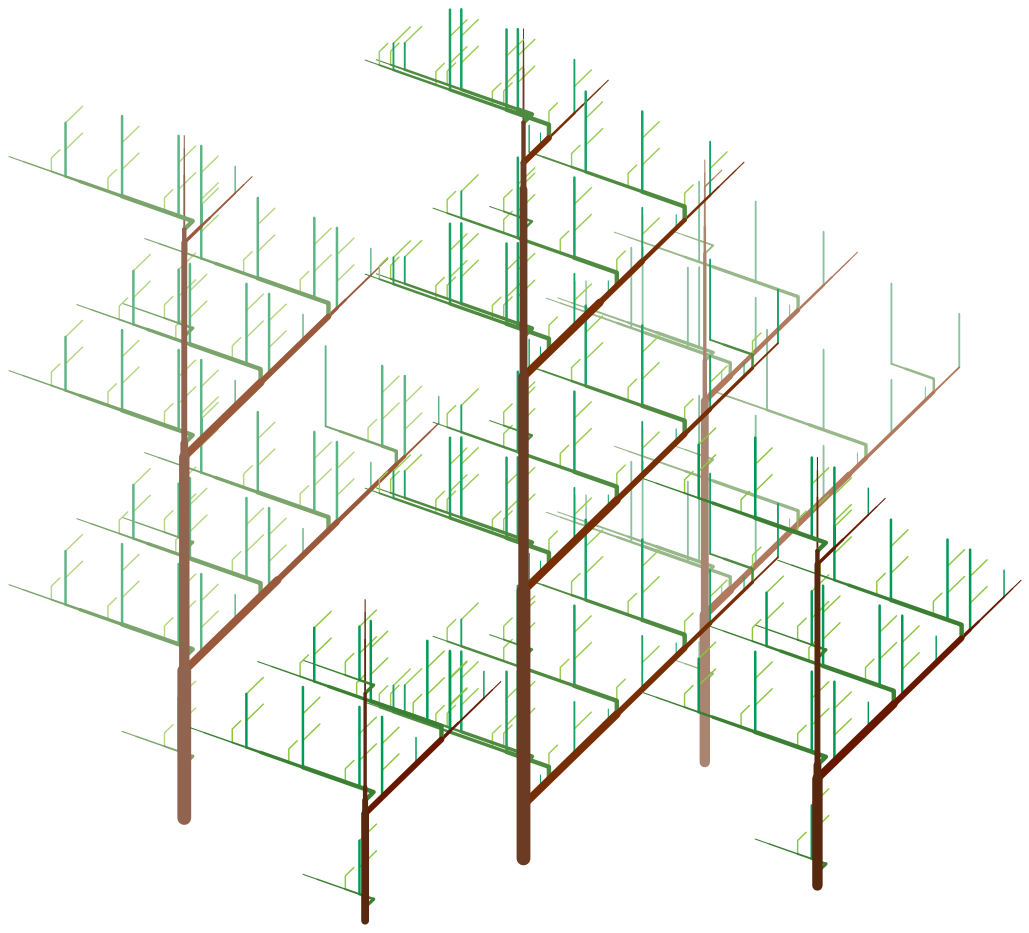

\centering

\caption{Part of the forest formed by the edges $\{x,p(x)\}$ in the proof of Theorem~\ref{t.counter}, drawn with every parent below its children.}
\label{f.forest}
\end{figure}

We first show how to change the conductances. The edges $\{x,p(x)\}$ form a forest (Figure~\ref{f.forest}). We cancel the error $6\Delta h(x)$ at each site $x$ by changing the conductance of the edge from $x$ to its parent. This also changes the harmonicity condition at the parent, so the errors accumulate up the forest: the change on the edge from $x$ to its parent must cancel the errors at $x$ and at all of its descendants. Lemma~\ref{l.correction} shows that these changes stay small when $|h|$ grows fast enough along the forest.

\begin{lemma}\label{l.correction}
Let $0<\varepsilon<\nf19$, and let $h\colon\Z^3\to\R$ satisfy $|h|\geq1$. For every $x$, let $p(x)$ be a neighbor of $x$ at which $|h|$ is largest, with ties broken by a fixed order of the six directions, and suppose that
\begin{equation}\label{e.growth}
|h(p(x))|\geq(3+2\sqrt2)|h(x)|\qquad\text{and}\qquad6|\Delta h(x)|\leq\varepsilon|h(p(x))|\qquad\text{for every }x\,.
\end{equation}
Then $h$ is $\bfa$-harmonic for some conductances $\bfa$ with $|\bfa-1|\leq9\varepsilon$. These are given by an explicit formula in $h$, and replacing $h$ by $h(\cdot+z)$ replaces $\bfa(x,y)$ by $\bfa(x+z,y+z)$.
\end{lemma}

\begin{proof}
Call $p(x)$ the parent of $x$, and $x$ a child of $p(x)$. We keep the conductance $1$ on every edge except the edges from sites to their parents. In Step~1, we define the correction $J(x)$ that the edge from $x$ to its parent must make: it cancels the error $6\Delta h(x)$ at $x$ and passes up the corrections from the children of $x$. In Step~2, we turn these corrections into conductances and check that $h$ is $\bfa$-harmonic.

\medskip\noindent\textbf{Step 1:} We prove that there is $J\colon\Z^3\to\R$ such that, for every $x$,
\begin{equation*}
J(x)=-6\Delta h(x)+\sum_{p(y)=x}J(y)\qquad\text{and}\qquad|J(x)|\leq\frac{\varepsilon|h(p(x))|}{10\sqrt2-14}\,.
\end{equation*}
Since $|h|\geq1$ grows by a factor at least $3+2\sqrt2$ from a child to its parent, each site has finitely many descendants and at most five children. So the recursion on the left in the above display defines $J$, starting from $J(x)=-6\Delta h(x)$ at the sites without children. By~\eqref{e.growth}, induction and the bound $|h(x)|\leq(3-2\sqrt2)|h(p(x))|$, we get
\begin{equation*}
|J(x)|\leq\varepsilon|h(p(x))|+\frac{5\varepsilon|h(x)|}{10\sqrt2-14}\leq\frac{\varepsilon|h(p(x))|}{10\sqrt2-14}\,.
\end{equation*}

\medskip\noindent\textbf{Step 2:} We conclude. We give each edge $\{x,p(x)\}$ the conductance
\begin{equation*}
\bfa(x,p(x))\coloneqq1+\frac{J(x)}{h(p(x))-h(x)}\,,
\end{equation*}
and every other edge conductance $1$. Since $|h(p(x))-h(x)|\geq(2\sqrt2-2)|h(p(x))|$, Step~1 gives $|\bfa(x,p(x))-1|\leq9\varepsilon$, and
\begin{equation*}
\sum_{|y-x|=1}\bfa(x,y)(h(y)-h(x))=6\Delta h(x)+J(x)-\sum_{p(y)=x}J(y)=0\,.\qedhere
\end{equation*}
\end{proof}

With Lemma~\ref{l.correction}, we can now prove Theorem~\ref{t.counter}.

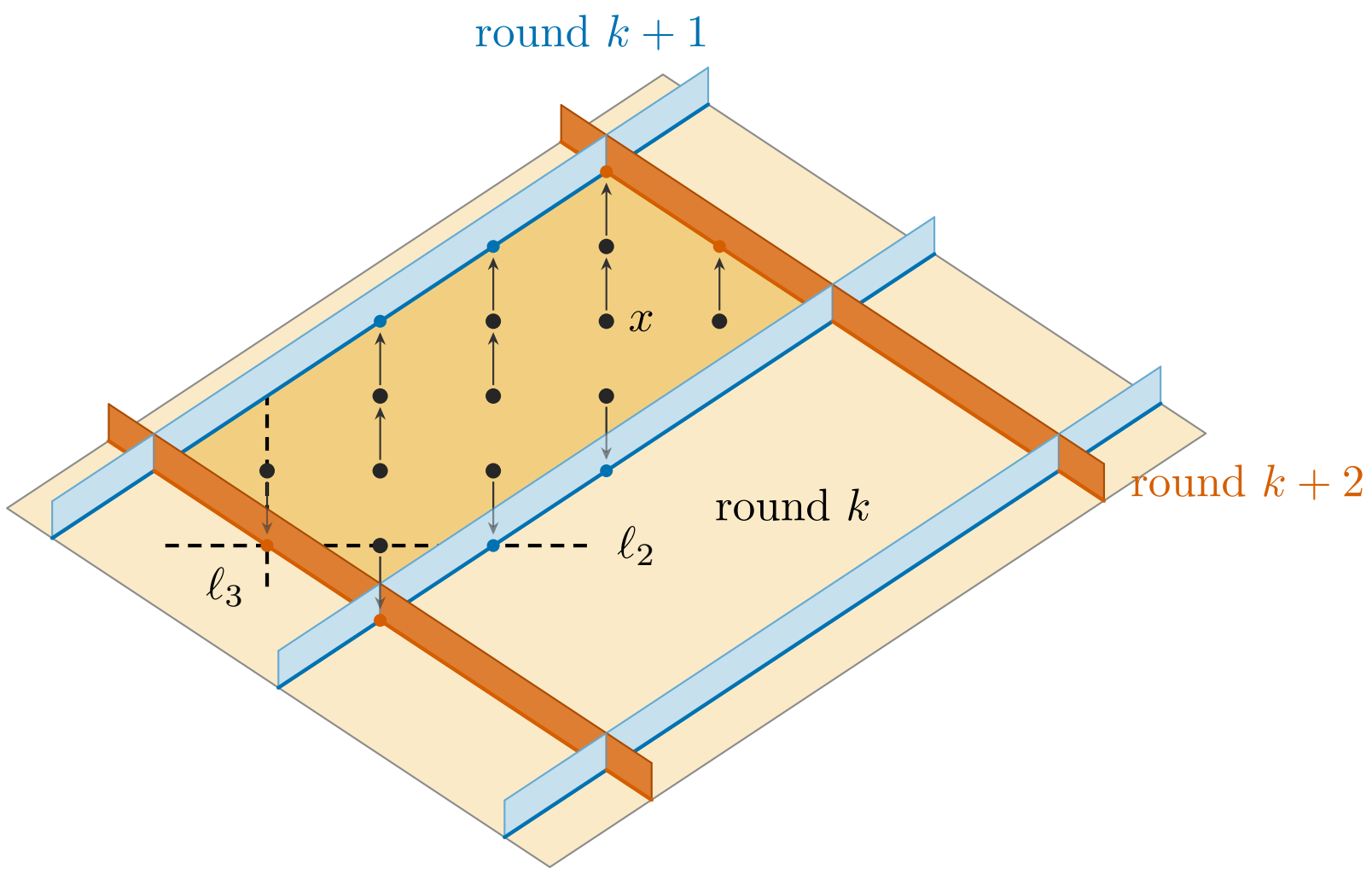
\begin{figure}[t]
\centering
\definecolor{OIorange}{HTML}{E69F00}\definecolor{OIblue}{HTML}{0072B2}\definecolor{OIverm}{HTML}{D55E00}
\begin{tikzpicture}[x={(0.0000cm,0.5794cm)},y={(0.4384cm,-0.1560cm)},z={(-0.4384cm,-0.1560cm)},
  line cap=round,line join=round,
  lab/.style={font=\small},
  planelab/.style={font=\small},
  blu/.style={OIblue},
  grn/.style={OIverm},
  site/.style={circle,fill=black!85,inner sep=1.2pt},
  ell/.style={densely dashed,thick,black,line cap=butt},
  dir/.style={-{Stealth[length=5pt]},thick},
  parent/.style={-{Stealth[length=2.8pt,width=2.4pt]},line width=0.4pt,black!80,line cap=butt,shorten >=2.4pt,shorten <=2.4pt},
  psite/.style={circle,inner sep=1pt}]
  \filldraw[fill=OIorange!22,draw=black!45,thin] (-3.3,1.5,-1.5) -- (1.5,-3.3,3.3) -- (7.3,2.5,-2.5) -- (2.5,7.3,-7.3) -- cycle;
  \fill[OIorange!50] (0,0,0) -- (4,4,-4) -- (6,2,-2) -- (2,-2,2) -- cycle;
  \draw[blu,thick] (-2.9,1.1,-1.1) -- (2.9,6.9,-6.9);
  \draw[blu,thick] (-0.9,-0.9,0.9) -- (4.9,4.9,-4.9);
  \draw[blu,thick] (1.1,-2.9,2.9) -- (6.9,2.9,-2.9);
  \draw[grn,thick] (-2.4,2.4,-2.4) -- (2.4,-2.4,2.4);
  \draw[grn,thick] (1.6,6.4,-6.4) -- (6.4,1.6,-1.6);
  \draw[ell] (1,-1.9,1.9) -- (1,1.9,-1.9);
  \draw[ell] (0.45,-1,1) -- (3,-1,1);
  \foreach \p in {(2,-1,1),(1,0,0),(2,0,0),(3,0,0),(2,1,-1),(3,1,-1),(4,1,-1),(3,2,-2),(4,2,-2),(5,2,-2),(4,3,-3)}
    \node[site] at \p {};
  \node[lab,right=1pt] at (4,2,-2) {$x$};
  \draw[parent] (2,-1,1) -- (1,-1,1);
  \draw[parent] (1,0,0) -- (0,0,0);
  \draw[parent] (2,0,0) -- (3,0,0);
  \draw[parent] (3,0,0) -- (4,0,0);
  \draw[parent] (2,1,-1) -- (1,1,-1);
  \draw[parent] (3,1,-1) -- (4,1,-1);
  \draw[parent] (4,1,-1) -- (5,1,-1);
  \draw[parent] (3,2,-2) -- (2,2,-2);
  \draw[parent] (4,2,-2) -- (5,2,-2);
  \draw[parent] (5,2,-2) -- (6,2,-2);
  \draw[parent] (4,3,-3) -- (5,3,-3);
  \filldraw[fill=OIverm!80,draw=OIverm!80,line width=0.2pt] (6,2,-2) -- (6.4,1.6,-1.6) -- (6.72,1.28,-1.92) -- (6.32,1.68,-2.32) -- cycle;
  \draw[OIverm!80!black,thin] (6.32,1.68,-2.32) -- (6.72,1.28,-1.92);
  \draw[black!45,thin] (6,2,-2) -- (6.32,1.68,-2.32);
  \draw[OIverm!80!black,thin] (6.4,1.6,-1.6) -- (6.72,1.28,-1.92);
  \draw[grn,thick] (6,2,-2) -- (6.4,1.6,-1.6);
  \filldraw[fill=OIverm!80,draw=OIverm!80,line width=0.2pt] (2,-2,2) -- (2.4,-2.4,2.4) -- (2.72,-2.72,2.08) -- (2.32,-2.32,1.68) -- cycle;
  \draw[OIverm!80!black,thin] (2.32,-2.32,1.68) -- (2.72,-2.72,2.08);
  \draw[black!45,thin] (2,-2,2) -- (2.32,-2.32,1.68);
  \draw[OIverm!80!black,thin] (2.4,-2.4,2.4) -- (2.72,-2.72,2.08);
  \draw[grn,thick] (2,-2,2) -- (2.4,-2.4,2.4);
  \filldraw[fill=OIblue!22,draw=OIblue!22,line width=0.2pt] (6,2,-2) -- (6.9,2.9,-2.9) -- (7.22,2.58,-3.22) -- (6.32,1.68,-2.32) -- cycle;
  \draw[OIblue!60,thin] (6.32,1.68,-2.32) -- (7.22,2.58,-3.22);
  \draw[black!45,thin] (6,2,-2) -- (6.32,1.68,-2.32);
  \draw[OIblue!60,thin] (6.9,2.9,-2.9) -- (7.22,2.58,-3.22);
  \draw[blu,thick] (6,2,-2) -- (6.9,2.9,-2.9);
  \filldraw[fill=OIblue!22,draw=OIblue!22,line width=0.2pt] (2,-2,2) -- (6,2,-2) -- (6.32,1.68,-2.32) -- (2.32,-2.32,1.68) -- cycle;
  \draw[OIblue!60,thin] (2.32,-2.32,1.68) -- (6.32,1.68,-2.32);
  \draw[black!45,thin] (2,-2,2) -- (2.32,-2.32,1.68);
  \draw[black!45,thin] (6,2,-2) -- (6.32,1.68,-2.32);
  \draw[blu,thick] (2,-2,2) -- (6,2,-2);
  \filldraw[fill=OIblue!22,draw=OIblue!22,line width=0.2pt] (1.1,-2.9,2.9) -- (2,-2,2) -- (2.32,-2.32,1.68) -- (1.42,-3.22,2.58) -- cycle;
  \draw[OIblue!60,thin] (1.42,-3.22,2.58) -- (2.32,-2.32,1.68);
  \draw[OIblue!60,thin] (1.1,-2.9,2.9) -- (1.42,-3.22,2.58);
  \draw[black!45,thin] (2,-2,2) -- (2.32,-2.32,1.68);
  \draw[blu,thick] (1.1,-2.9,2.9) -- (2,-2,2);
  \filldraw[fill=OIverm!80,draw=OIverm!80,line width=0.2pt] (4,4,-4) -- (6,2,-2) -- (6.32,1.68,-2.32) -- (4.32,3.68,-4.32) -- cycle;
  \draw[OIverm!80!black,thin] (4.32,3.68,-4.32) -- (6.32,1.68,-2.32);
  \draw[black!45,thin] (4,4,-4) -- (4.32,3.68,-4.32);
  \draw[black!45,thin] (6,2,-2) -- (6.32,1.68,-2.32);
  \draw[grn,thick] (4,4,-4) -- (6,2,-2);
  \filldraw[fill=OIverm!80,draw=OIverm!80,line width=0.2pt] (0,0,0) -- (2,-2,2) -- (2.32,-2.32,1.68) -- (0.32,-0.32,-0.32) -- cycle;
  \draw[OIverm!80!black,thin] (0.32,-0.32,-0.32) -- (2.32,-2.32,1.68);
  \draw[black!45,thin] (0,0,0) -- (0.32,-0.32,-0.32);
  \draw[black!45,thin] (2,-2,2) -- (2.32,-2.32,1.68);
  \draw[grn,thick] (0,0,0) -- (2,-2,2);
  \filldraw[fill=OIblue!22,draw=OIblue!22,line width=0.2pt] (4,4,-4) -- (4.9,4.9,-4.9) -- (5.22,4.58,-5.22) -- (4.32,3.68,-4.32) -- cycle;
  \draw[OIblue!60,thin] (4.32,3.68,-4.32) -- (5.22,4.58,-5.22);
  \draw[black!45,thin] (4,4,-4) -- (4.32,3.68,-4.32);
  \draw[OIblue!60,thin] (4.9,4.9,-4.9) -- (5.22,4.58,-5.22);
  \draw[blu,thick] (4,4,-4) -- (4.9,4.9,-4.9);
  \filldraw[fill=OIblue!22,draw=OIblue!22,line width=0.2pt] (0,0,0) -- (4,4,-4) -- (4.32,3.68,-4.32) -- (0.32,-0.32,-0.32) -- cycle;
  \draw[OIblue!60,thin] (0.32,-0.32,-0.32) -- (4.32,3.68,-4.32);
  \draw[black!45,thin] (0,0,0) -- (0.32,-0.32,-0.32);
  \draw[black!45,thin] (4,4,-4) -- (4.32,3.68,-4.32);
  \draw[blu,thick] (0,0,0) -- (4,4,-4);
  \filldraw[fill=OIblue!22,draw=OIblue!22,line width=0.2pt] (-0.9,-0.9,0.9) -- (0,0,0) -- (0.32,-0.32,-0.32) -- (-0.58,-1.22,0.58) -- cycle;
  \draw[OIblue!60,thin] (-0.58,-1.22,0.58) -- (0.32,-0.32,-0.32);
  \draw[OIblue!60,thin] (-0.9,-0.9,0.9) -- (-0.58,-1.22,0.58);
  \draw[black!45,thin] (0,0,0) -- (0.32,-0.32,-0.32);
  \draw[blu,thick] (-0.9,-0.9,0.9) -- (0,0,0);
  \filldraw[fill=OIverm!80,draw=OIverm!80,line width=0.2pt] (2,6,-6) -- (4,4,-4) -- (4.32,3.68,-4.32) -- (2.32,5.68,-6.32) -- cycle;
  \draw[OIverm!80!black,thin] (2.32,5.68,-6.32) -- (4.32,3.68,-4.32);
  \draw[black!45,thin] (2,6,-6) -- (2.32,5.68,-6.32);
  \draw[black!45,thin] (4,4,-4) -- (4.32,3.68,-4.32);
  \draw[grn,thick] (2,6,-6) -- (4,4,-4);
  \filldraw[fill=OIverm!80,draw=OIverm!80,line width=0.2pt] (-2,2,-2) -- (0,0,0) -- (0.32,-0.32,-0.32) -- (-1.68,1.68,-2.32) -- cycle;
  \draw[OIverm!80!black,thin] (-1.68,1.68,-2.32) -- (0.32,-0.32,-0.32);
  \draw[black!45,thin] (-2,2,-2) -- (-1.68,1.68,-2.32);
  \draw[black!45,thin] (0,0,0) -- (0.32,-0.32,-0.32);
  \draw[grn,thick] (-2,2,-2) -- (0,0,0);
  \filldraw[fill=OIblue!22,draw=OIblue!22,line width=0.2pt] (2,6,-6) -- (2.9,6.9,-6.9) -- (3.22,6.58,-7.22) -- (2.32,5.68,-6.32) -- cycle;
  \draw[OIblue!60,thin] (2.32,5.68,-6.32) -- (3.22,6.58,-7.22);
  \draw[black!45,thin] (2,6,-6) -- (2.32,5.68,-6.32);
  \draw[OIblue!60,thin] (2.9,6.9,-6.9) -- (3.22,6.58,-7.22);
  \draw[blu,thick] (2,6,-6) -- (2.9,6.9,-6.9);
  \filldraw[fill=OIblue!22,draw=OIblue!22,line width=0.2pt] (-2,2,-2) -- (2,6,-6) -- (2.32,5.68,-6.32) -- (-1.68,1.68,-2.32) -- cycle;
  \draw[OIblue!60,thin] (-1.68,1.68,-2.32) -- (2.32,5.68,-6.32);
  \draw[black!45,thin] (-2,2,-2) -- (-1.68,1.68,-2.32);
  \draw[black!45,thin] (2,6,-6) -- (2.32,5.68,-6.32);
  \draw[blu,thick] (-2,2,-2) -- (2,6,-6);
  \filldraw[fill=OIblue!22,draw=OIblue!22,line width=0.2pt] (-2.9,1.1,-1.1) -- (-2,2,-2) -- (-1.68,1.68,-2.32) -- (-2.58,0.78,-1.42) -- cycle;
  \draw[OIblue!60,thin] (-2.58,0.78,-1.42) -- (-1.68,1.68,-2.32);
  \draw[OIblue!60,thin] (-2.9,1.1,-1.1) -- (-2.58,0.78,-1.42);
  \draw[black!45,thin] (-2,2,-2) -- (-1.68,1.68,-2.32);
  \draw[blu,thick] (-2.9,1.1,-1.1) -- (-2,2,-2);
  \filldraw[fill=OIverm!80,draw=OIverm!80,line width=0.2pt] (1.6,6.4,-6.4) -- (2,6,-6) -- (2.32,5.68,-6.32) -- (1.92,6.08,-6.72) -- cycle;
  \draw[OIverm!80!black,thin] (1.92,6.08,-6.72) -- (2.32,5.68,-6.32);
  \draw[OIverm!80!black,thin] (1.6,6.4,-6.4) -- (1.92,6.08,-6.72);
  \draw[black!45,thin] (2,6,-6) -- (2.32,5.68,-6.32);
  \draw[grn,thick] (1.6,6.4,-6.4) -- (2,6,-6);
  \filldraw[fill=OIverm!80,draw=OIverm!80,line width=0.2pt] (-2.4,2.4,-2.4) -- (-2,2,-2) -- (-1.68,1.68,-2.32) -- (-2.08,2.08,-2.72) -- cycle;
  \draw[OIverm!80!black,thin] (-2.08,2.08,-2.72) -- (-1.68,1.68,-2.32);
  \draw[OIverm!80!black,thin] (-2.4,2.4,-2.4) -- (-2.08,2.08,-2.72);
  \draw[black!45,thin] (-2,2,-2) -- (-1.68,1.68,-2.32);
  \draw[grn,thick] (-2.4,2.4,-2.4) -- (-2,2,-2);
  \draw[parent,opacity=0.55] (2,-1,1) -- (1,-1,1);
  \draw[parent,opacity=0.55] (1,0,0) -- (0,0,0);
  \draw[parent,opacity=0.55] (2,1,-1) -- (1,1,-1);
  \draw[parent,opacity=0.55] (3,2,-2) -- (2,2,-2);
  \node[psite,fill=OIverm] at (0,0,0) {};
  \node[psite,fill=OIverm] at (1,-1,1) {};
  \node[psite,fill=OIblue] at (1,1,-1) {};
  \node[psite,fill=OIblue] at (2,2,-2) {};
  \node[psite,fill=OIblue] at (4,0,0) {};
  \node[psite,fill=OIblue] at (5,1,-1) {};
  \node[psite,fill=OIverm] at (5,3,-3) {};
  \node[psite,fill=OIverm] at (6,2,-2) {};
  \node[planelab,anchor=center] at (1.55,3.65,-3.65) {round $k$};
  \node[planelab,blu,anchor=south east,xshift=4pt] at (7.22,2.58,-3.22) {round $k+1$};
  \node[planelab,grn,anchor=west,xshift=2pt] at (1.76,6.24,-6.56) {round $k+2$};
  \node[lab,right=1pt] at (1,1.9,-1.9) {$\ell_2$};
  \node[lab,left=1pt] at (0.45,-1,1) {$\ell_3$};
\end{tikzpicture}
\caption{The planes marked in rounds $k+1$ (blue) and $k+2$ (dark orange) cut a plane marked in round $k$ (orange) into parallelograms. In the parallelogram containing $x$ (shaded), the dashed lines mark the smallest second and third coordinates $\ell_2$ and $\ell_3$ of its sites (here $k=1$ and $i=2$). Arrows point from each site to its parent.}
\label{f.planes}
\end{figure}

\begin{proof}[Proof of Theorem~\ref{t.counter}]
Let $\lambda>1$. It suffices to treat $d=3$, since an example on $\Z^3$ extends to $\Z^d$ by stacking copies. We glue together the harmonic functions
\begin{equation}\label{e.planefns}
x\mapsto(-1)^{x_{i+1}}(3+2\sqrt2)^{x_i}\,\mathbf 1\{x_{i+1}+x_{i+2}=s\}\qquad(i\in\{1,2,3\},\ s\in\Z)\,,
\end{equation}
with indices taken modulo $3$, in a stationary way so that~\eqref{e.growth} holds with $9\varepsilon<1-\lambda^{-1}$, and then apply Lemma~\ref{l.correction}. Each of these functions is supported on a diagonal plane, and for each $i$ these planes form a family of parallel planes.

To make the gluing stationary, we construct $h$ as follows. Let $U_0,U_1,\dots$ be independent random variables, with $U_k$ uniform on $\{0,\dots,2^k-1\}$, and fix $\beta >108/(1-\lambda^{-1})$.
\begin{enumerate}
\item[(a)] For $k\geq0$, let $i\in\{1,2,3\}$ with $i\equiv k+1$ modulo $3$. In round $k$, we mark the planes $\{x_{i+1}+x_{i+2}=s\}$ with $s\equiv U_k$ modulo $2^k$. A site is marked in round $k$ if it lies on one of these planes. In round $0$, every plane $\{x_2+x_3=s\}$ is marked. A site is marked in round $k$ with probability $2^{-k}$, so by the Borel--Cantelli lemma, almost surely every site is marked in only finitely many rounds.
\item[(b)] Let $x$ be a site, let $k$ be the last round in which $x$ is marked, and let $i\in\{1,2,3\}$ with $i\equiv k+1$ modulo $3$. The planes marked in rounds $k+1$ and $k+2$ cut the plane through $x$ marked in round $k$ into parallelograms (Figure~\ref{f.planes}). Let $\pgram{x}$ be the set of sites strictly inside the parallelogram containing $x$.
\item[(c)] Let $\ell_i$ and $\ell_{i+1}$ be the smallest $i$th and $(i+1)$st coordinates of the sites in $\pgram{x}$, and set
\begin{equation}\label{e.hdef}
h(x)\coloneqq(-1)^{x_{i+1}-\ell_{i+1}}\beta^k(3+2\sqrt2)^{3(2^k-1)+x_i-\ell_i}\,.
\end{equation}
\end{enumerate}
Translating the lattice shifts the random residue classes without changing their law, so $h$ is stationary.

\medskip\noindent\textbf{Step 1:} We prove that $|h|\geq1$ and that if $y$ was last marked in a later round than $x$, then
\begin{equation}\label{e.separation}
|h(x)|<\beta^{-1}|h(y)|\,.
\end{equation}
If $x$ was last marked in round $k$, then it lies on none of the planes marked in rounds $k+1$ and $k+2$, so $\pgram{x}$ is well defined. The set $\pgram{x}$ lies between consecutive planes marked in rounds $k+1$ and $k+2$, which are $2^{k+1}$ and $2^{k+2}$ planes apart, so $0\leq x_i-\ell_i<3\cdot2^k$. Hence $\beta^k(3+2\sqrt2)^{3(2^k-1)}\leq|h(x)|<\beta^{-1}\cdot \beta^{k+1}(3+2\sqrt2)^{3(2^{k+1}-1)}$, which gives $|h|\geq1$ and~\eqref{e.separation}.

\medskip\noindent\textbf{Step 2:} Near each site, $h$ is dominated by its values on the planes marked in the latest round, and on these planes $h$ agrees with a sum of the harmonic functions~\eqref{e.planefns}. We make this precise. Fix a site $x$, let $\mathcal N(x)\coloneqq\{x\}\cup\{x\pm e_j:j=1,2,3\}$ be the set consisting of $x$ and its six neighbors, let $k$ be the last round in which a site of $\mathcal N(x)$ is marked, and let $\mathcal N^{(k)}(x)$ be the set of sites of $\mathcal N(x)$ last marked in round $k$. Let $i\in\{1,2,3\}$ with $i\equiv k+1$ modulo $3$. We construct a harmonic function $u\colon\Z^3\to\R$ that equals $h$ on $\mathcal N^{(k)}(x)$ and vanishes on $\mathcal N(x)\setminus\mathcal N^{(k)}(x)$.

By construction, the sites of $\mathcal N^{(k)}(x)$ lie on parallel planes marked in round $k$, and those on each such plane $\Pi$ lie in a single parallelogram. Let $\ell_i^\Pi$ and $\ell_{i+1}^\Pi$ be the values of $\ell_i$ and $\ell_{i+1}$ for this parallelogram, and set
\begin{equation*}
u(y)\coloneqq\sum_\Pi (-1)^{y_{i+1}-\ell_{i+1}^\Pi}\beta^k(3+2\sqrt2)^{3(2^k-1)+y_i-\ell_i^\Pi}\,\mathbf 1\{y\in \Pi\}\,,
\end{equation*}
summing over the planes $\Pi$ marked in round $k$ that meet $\mathcal N(x)$. Each term is a constant multiple of one of the harmonic functions~\eqref{e.planefns}, so $\Delta u=0$. By the definition~\eqref{e.hdef}, $u=h$ on $\mathcal N^{(k)}(x)$, and $u=0$ on $\mathcal N(x)\setminus\mathcal N^{(k)}(x)$, since these sites lie on none of these planes.

\medskip\noindent\textbf{Step 3:} We prove that~\eqref{e.growth} holds with $\varepsilon\coloneqq12/\beta$, and conclude. If $x\in\mathcal N^{(k)}(x)$, then $x+e_i\in\pgram{x}$, and $|h(x+e_i)|=(3+2\sqrt2)|h(x)|$. Otherwise, some neighbor of $x$ lies in $\mathcal N^{(k)}(x)$, and $|h|>\beta|h(x)|$ there by~\eqref{e.separation}. In both cases, $|h(p(x))|\geq(3+2\sqrt2)|h(x)|$ and $p(x)\in\mathcal N^{(k)}(x)$ (Figure~\ref{f.planes}). Since $\Delta u(x)=0$ and $|h-u|<\beta^{-1}|h(p(x))|$ on $\mathcal N(x)$ by~\eqref{e.separation}, we get $6|\Delta h(x)|\leq12\beta^{-1}|h(p(x))|$. Since $9\varepsilon=108/\beta <1-\lambda^{-1}$, Lemma~\ref{l.correction} now completes the proof.
\end{proof}

\appendix

\section{Large-scale analyticity}\label{s.analytic}

We recall that harmonic functions on $\Z^d$ satisfy large-scale analyticity, in the sense of~\citet{AKS}. Related quantitative estimates for discrete harmonic functions on $\Z^d$ include the three balls theorem of~\citet{GuadieMalinnikova} and the propagation of smallness of~\citet{LippnerMangoubi}.

\begin{lemma}\label{l.harmonic}
For every $\alpha\in(0,1)$ there are $c>0$ and $N_0$, depending only on $d$ and $\alpha$, with the following property. Let $N\geq N_0$, let $1\leq m\leq cN$ be an integer, and let $u\colon\Z^d\to\R$ be harmonic in $Q_{3N}$. There is a real polynomial $P\colon\R^d\to\R$ of total degree at most $m$ such that
\begin{equation*}
\|u-P\|_{L^\infty(Q_{cN})}\leq\alpha^m\|u\|_{\underline L^2(Q_{3N})}\,.
\end{equation*}
\end{lemma}

\begin{proof}
Since the proof is nearly identical to that of~\citet[Lemma~A.2]{BCG}, we only sketch the argument. For a multi-index $\nu$, let $D^\nu$ denote the mixed forward difference of order $\nu_j$ in each direction $e_j$. Iterating the discrete Caccioppoli inequality \citep[Proposition~12]{BDCKY} and applying a discrete Moser estimate \citep[Proposition~5.3]{Delmotte} to the subharmonic function $|D^\nu u|$,\footnote{\citet{Delmotte} writes $Lg(x)=\sum_yC_{xy}(g(x)-g(y))$. The hypothesis $Lg\geq0$ of~\citet[Proposition~5.3]{Delmotte} is misprinted: it should be $Lg\leq0$, which says that $g$ is subharmonic, as in its proof.} we obtain, for $|\nu|=s$ with $1\leq s\leq N/C$,
\begin{equation}\label{e.mixed}
\|D^\nu u\|_{L^\infty(Q_N)}\leq CN^{-\nf d2}\|D^\nu u\|_{L^2(Q_{2N})}\leq\left(\frac{Cs}{N}\right)^s\|u\|_{\underline L^2(Q_{3N})}\,.
\end{equation}
Here we use the fact that the differences of a harmonic function are harmonic. Let $P$ be the polynomial of degree at most $m$ with $D^\nu P(0)=D^\nu u(0)$ for $|\nu|\leq m$. For $x\in Q_{cN}$, summing $m+1$ times along a shortest coordinate path from $0$ to $x$ yields
\begin{equation*}
|u(x)-P(x)|\leq\binom{|x|_1+m}{m+1}\max_{|\nu|=m+1}\|D^\nu u\|_{L^\infty(Q_{cN})}\,.
\end{equation*}
The lemma follows by applying~\eqref{e.mixed} with $s=m+1$, choosing $c$ sufficiently small and $N_0$ sufficiently large.
\end{proof}

\renewcommand{\bibfont}{\small}
\setlength{\bibsep}{1pt}

\bibliographystyle{plainnat}
\bibliography{references}

\end{document}